\documentclass{slides}

\usepackage{amsmath,amssymb}
\usepackage{latexsym}

\newcommand{\zz}{\mathbb{Z}}

\newcommand{\bH}{\mathbf{H} }
\newcommand{\bHnm}{\mathbf{H}_{n,m}}

\newcommand{\ba}{\backslash }
\newcommand{\GL}{GL_n(\mathbb{R})}

\newcommand\Rg{{\mathcal R}_n}

\begin{document}
\title{\Huge \bf Remark on Harmonic Analysis on Siegel-Jacobi Space}
\author{{\Large \bf Jae-Hyun Yang}\\ \vskip 0.5cm {\large \bf Inha University} \\
jhyang@inha.ac.kr}
\date{Department of Mathematics\\ Kyoto University\\
 Kyoto, Japan
\\ June 23 (Tue), 2009}

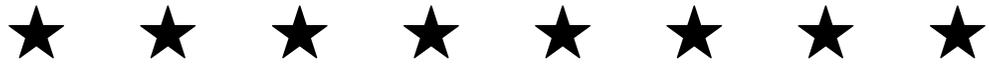
\def\lrt{\longrightarrow}
\def\CP{{\bf{\cal P}}_n}
\def\CPR{{\cal P}_n\times {\Bbb R}^{(m,n)}}
\def\SP{S{\cal P}_n}
\def\GL{GL_n(\Bbb R)}
\def\Rmn{{\Bbb R}^{(m,n)}}
\def\BR{\Bbb R}
\def\BC{\Bbb C}
\newcommand\BZ{\mathbb Z}
\def\g{\gamma}
\def\G{\Gamma}
\def\O{\Omega}
\def\o{\omega}
\def\la{\lambda}
\def\BD{\Bbb D}
\def\BH{\Bbb H}
\def\Cmn{{\Bbb C}^{(m,n)}}
\def\PZ{ {{\partial}\over {\partial Z}} }
\def\PW{ {{\partial}\over {\partial W}} }
\def\PZB{ {{\partial}\over {\partial{\overline Z}}} }
\def\PWB{ {{\partial}\over {\partial{\overline W}}} }
\def\PO{ {{\partial}\over {\partial \Omega}} }
\def\PE{ {{\partial}\over {\partial \eta}} }
\def\POB{ {{\partial}\over {\partial{\overline \Omega}}} }
\def\PEB{ {{\partial}\over {\partial{\overline \eta}}} }
\def\PX{ {{\partial}\over{\partial X}} }
\def\PY{ {{\partial}\over {\partial Y}} }
\def\PU{ {{\partial}\over{\partial U}} }
\def\PV{ {{\partial}\over{\partial V}} }
\def\HC{\mathbf{H}_n\times \BC^{(m,n)} }

\newcommand\tr{\textrm{tr}}
\newcommand\Om{\Omega}
\newcommand\bz{d{\overline Z}}
\newcommand\bo{d{\overline \Omega}}
\newcommand\om{\omega}
\newcommand\Dnm{{\mathbb D}_{n,m}}
\newcommand\Hn{{\mathbb H}_n}
\newcommand\ka{\kappa}
\newcommand\w{\wedge}

\newcommand\POBS{ {{\partial\ \,}\over {\partial{\overline \Omega}_*} } }
\newcommand\PZBS{ {{\partial\ \,}\over {\partial{\overline Z_*}}} }
\newcommand\PXS{ {{\partial\ \,}\over{\partial X_*}} }
\newcommand\PYS{ {{\partial\ \,}\over {\partial Y_*}} }
\newcommand\PUS{ {{\partial\ \,}\over{\partial U_*}} }
\newcommand\PVS{ {{\partial\ \,}\over{\partial V_*}} }
\newcommand\POS{ {{\partial\ \,}\over{\partial \Omega_*}} }
\newcommand\PZS{ {{\partial\ \,}\over{\partial Z_*}} }

\newcommand\Dg{{\mathbb D}_n}
\newcommand\Hg{{\mathbb H}_n}
\newcommand\bDn{{\mathbf D}_n}
\newcommand\bDnm{{\mathbf D}_{n,m}}

\newcommand\OW{\overline{W}}
\newcommand\OP{\overline{P}}
\newcommand\OQ{\overline{Q}}
\newcommand\OVW{\overline W}

\newcommand\ot{\overline\eta}

\newcommand\Fgh{{\mathcal F}_{n,m}}

\maketitle

%%%%%%%%%%%%%%%%%%%%%%%%%%%%%%%%%%%%%%%%%%%%%%%%%%%%%%%%%%%%%%%%%%%%%%%%%%%%%%%%%%%%%%%%%%%%%%%%%%%%%%%%%%%%%%%%%%%%%%%%%%%%%%%%%%%%
%%%%%%%%%%%%%%%%%%%%%%%%%%%%%%%%%%%%%%%%%%%%%%%%%%%%%%%%%%%%%%%%%%%%%%%%%%%%%%%%%%%%%%%%%%%%%%%%%%%%%%%%%%%%%%%%%%%%%%%%%%%%%%%%%%%%
%\newpage
%\begin{slide}
%\begin{center}
%{\large \bf 1. Langlands Functoriality Conjecture (briefly LFC)}
%\end{center}

\newpage
\begin{slide}
\begin{center}
{\Large $ \textbf{$\bigstar\ \ \bigstar\ \ \bigstar\ \ \bigstar\ \ \bigstar
\ \ \bigstar\ \ \bigstar\ \ \bigstar$}$}
\end{center}

\indent\ \ \ My work was inspired by the spirit of the great
number theorists of the 20th century

\vskip 1cm $ \textbf{Carl Ludiwig Siegel (1896-1981)}$
\vskip 0.1cm
$ \textbf{Andr{\'e} Weil (1906-1998)}$
\vskip 0.1cm
$ \textbf{Hans Maass (1911-1992)}$
\vskip 0.1cm $ \textbf{Atle Selberg (1917-2007)}$
\vskip 0.1cm $ \textbf{Robert P. Langlands (1936-\ )}$

\vskip 1cm

************************************\\
************************************\\
************************************

\end{slide}

%%%%%%%%%%%%%%%%%%%%%%%%%%%%%%%%%%%%%%%%%%%%%%%%%%%%%%%%%%%%%%%%%%%%%%%%%%%%%%%%%%%%%%%%%%%%%%%%%%%%%%%%%%%%%%%%%%%%%%%%%%%%%%%%%%%%%%%
%%%%%%%%%%%%%%%%%%%%%%%%%%%%%%%%%%%%%%%%%%%%%%%%%%%%%%%%%%%%%%%%%%%%%%%%%%%%%%%%%%%%%%%%%%%%%%%%%%%%%%%%%%%%%%%%%%%%%%%%%%%%%%%%%%%%%%%
%%%%%%%%%%%%%%%%%%%%%%%%%%%%%%%%%%%%%%%%%%%%%%%%%%%%%%%%%%%%%%%%%%%%%%%%%%%%%%%%%%%%%%%%%%%%%%%%%%%%%%%%%%%%%%%%%%%%%%%%%%%%%%%%%%%%%%%

\newpage
\begin{slide}

[A] $ \textbf{C. L. Siegel}$, {\em Symplectic Geometry,} Amer. J.
Math. {\bf 65} (1943), 1-86; Academic Press, New York and London
(1964); Gesammelte Abhandlungen, no.\,\,{\bf 41,\ vol.\,II},
Springer-Verlag (1966), 274-359.

[B] $ \textbf{H. Maass}$, {\em {\"U}ber eine neue Art von
nichtanalytischen automorphen Funktionen und die Bestimmung
Dirichletscher Reihen durch Funtionalgleichungen,} Math. Ann. $
\textbf{121}$ (1949), 141-183.

[C] $ \textbf{A. Selberg}$,  {\em Harmonic analysis and
discontinuous groups in weakly symmetric Riemannian spaces with
applications to Dirichlet series,} J. Indian Math. Soc. B.  {\bf
20} (1956), 47-87.

[D] $ \textbf{A. Weil}$, {\em Sur certains groupes d'operateurs
unitaires (French)}, Acta Math. $ \textbf{111}$ (1964), 143--211.

\end{slide}

%%%%%%%%%%%%%%%%%%%%%%%%%%%%%%%%%%%%%%%%%%%%%%%%%%%%%%%%%%%%%%%%%%%%%%%%%%%%%%%%%%%%%%%%%%%%%%%%%%%%%%%%%%%%%%%%%%%%%%%%%%%%%%%%%%%%%%%
%%%%%%%%%%%%%%%%%%%%%%%%%%%%%%%%%%%%%%%%%%%%%%%%%%%%%%%%%%%%%%%%%%%%%%%%%%%%%%%%%%%%%%%%%%%%%%%%%%%%%%%%%%%%%%%%%%%%%%%%%%%%%%%%%%%%%%%
%%%%%%%%%%%%%%%%%%%%%%%%%%%%%%%%%%%%%%%%%%%%%%%%%%%%%%%%%%%%%%%%%%%%%%%%%%%%%%%%%%%%%%%%%%%%%%%%%%%%%%%%%%%%%%%%%%%%%%%%%%%%%%%%%%%%%%%

\newpage
\begin{slide}
\begin{center}
{\large \bf Introduction}
\end{center}

Let
$$\bH_n=\left\{\,\O\in \BC^{(n,n)}\,|\ \O=\,{}^t\O,\ \mathrm {Im}\,\O>0\,\right\}$$
be the Siegel upper half plane and let
\begin{equation*}
\bHnm=\bH_n\times \BC^{(m,n)} \end{equation*} \noindent be the
Siegel-Jacobi space.

{\bf Notations\,:} Here $F^{(m,n)}$ denotes the set of all
$m\times n$ matrices with entries in a commutative ring $F$ and
${}^tA$ denotes the transpose of a matrix $A$. For an $n\times m$
matrix $B$ and an $n\times n$ matrix $A$, we write $A[B]=\,{}^tBAB.$

Let
\begin{equation*}
Sp(n,\BR)=\left\{\, M\in \BR^{(2n,2n)}\,|\ {}^tMJ_nM=J_n\ \right\}
\end{equation*}

\noindent be the symplectic group of degree $n$, where
\begin{equation*}
J_n=\begin{pmatrix} \ 0 & I_n \\ -I_n & 0 \end{pmatrix}.
\end{equation*}

\end{slide}

%%%%%%%%%%%%%%%%%%%%%%%%%%%%%%%%%%%%%%%%%%%%%%%%%%%%%%%%%%%%%%%%%%%%%%%%%%%%%%%%%%%%%%%%%%%%%%%%%%%%%%%%%%%%%%%%%%%%%%%%%%%%%%%%%%%%%%%
%%%%%%%%%%%%%%%%%%%%%%%%%%%%%%%%%%%%%%%%%%%%%%%%%%%%%%%%%%%%%%%%%%%%%%%%%%%%%%%%%%%%%%%%%%%%%%%%%%%%%%%%%%%%%%%%%%%%%%%%%%%%%%%%%%%%%%%
%%%%%%%%%%%%%%%%%%%%%%%%%%%%%%%%%%%%%%%%%%%%%%%%%%%%%%%%%%%%%%%%%%%%%%%%%%%%%%%%%%%%%%%%%%%%%%%%%%%%%%%%%%%%%%%%%%%%%%%%%%%%%%%%%%%%%%%

\begin{slide}
Then $Sp(n,\BR)$ acts on $\bH_n$ transitively by
\begin{align}
M\circ \O=(A\O+B)(C\O+D)^{-1},
\end{align}
where $M=\begin{pmatrix} A & B\\ C & D \end{pmatrix}\in Sp(n,\BR)$
and $\O\in\bH_n.$ Therefore
$$Sp(n,\BR)/U(n)\cong \bH_n$$
is a (Hermitian) symmetric space.
\par \ \ \ Let
$$H_{\BR}^{(n,m)}=\big\{ (\lambda,\mu,\kappa)\,\big| \
\lambda,\mu\in \BR^{(m,n)},\ \kappa\in \BR^{(m,m)}\,\big\}$$ be
the Heisenberg group. Let
$$G^J=Sp(n,\BR)\ltimes H_{\BR}^{(n,m)}$$
be the ${\mathbf{Jacobi\ group}}$ with the multiplication law
\begin{align*} &\ \ \ (M_0,(\lambda_0,\mu_0,\kappa_0))\cdot (M,(\lambda,\mu,\kappa))\\
& =\Big(M_0M,\big(\tilde{\lambda}_0+\lambda,\tilde{\mu}_0+\mu,
\kappa_0+\kappa+\tilde{\lambda}_0{}^t\!\mu-\tilde{\mu}_0{}^t\!\lambda\big)\Big),\end{align*}
\end{slide}

%%%%%%%%%%%%%%%%%%%%%%%%%%%%%%%%%%%%%%%%%%%%%%%%%%%%%%%%%%%%%%%%%%%%%%%%%%%%%%%%%%%%%%%%%%%%%%%%%%%%%%%%%%%%%%%%%%%%%%%%%%%%%%%%%%%%%%%
%%%%%%%%%%%%%%%%%%%%%%%%%%%%%%%%%%%%%%%%%%%%%%%%%%%%%%%%%%%%%%%%%%%%%%%%%%%%%%%%%%%%%%%%%%%%%%%%%%%%%%%%%%%%%%%%%%%%%%%%%%%%%%%%%%%%%%%
%%%%%%%%%%%%%%%%%%%%%%%%%%%%%%%%%%%%%%%%%%%%%%%%%%%%%%%%%%%%%%%%%%%%%%%%%%%%%%%%%%%%%%%%%%%%%%%%%%%%%%%%%%%%%%%%%%%%%%%%%%%%%%%%%%%%%%%
\begin{slide}

where $(\tilde{\lambda}_0,\tilde{\mu}_0)=(\lambda_0,\mu_0)M.$ Then
$G^J$ acts on the {\bf Siegel-Jacobi space} $\bHnm$ transitively
by
\begin{align} &\ \big(M,(\lambda,\mu,\kappa)\big)\cdot (\O,Z)\\
& =\big(M\circ \O, (Z+\lambda
\O+\mu)(C\O+D)^{-1}\big),\notag\end{align} where
$M=\begin{pmatrix} A & B\\ C & D \end{pmatrix}\in Sp(n,\BR),\
(\lambda,\mu,\kappa)\in H_{\BR}^{(n,m)}$ and $(\O,Z)\in \bHnm.$
Thus
$$G^J/K^J\cong \bHnm $$
is a {\bf non-reductive} complex manifold, where
$$K^J=U(n)\times \textrm{Sym}(n,\BR).$$
Let $\G_*$ be an arithmetic subgroup of $Sp(n,\BR)$ and
$\G_*^J=\G_*\ltimes H_{\zz}^{(n,m)}$. For instance,
$\G_*=Sp(n,\zz).$ Here
$$H_{\zz}^{(n,m)}=\big\{ (\lambda,\mu,\kappa)\in
H_{\BR}^{(n,m)}\,\big|\ \lambda,\mu,\kappa\ \textrm{integral}\
\big\}.$$

\end{slide}

%%%%%%%%%%%%%%%%%%%%%%%%%%%%%%%%%%%%%%%%%%%%%%%%%%%%%%%%%%%%%%%%%%%%%%%%%%%%%%%%%%%%%%%%%%%%%%%%%%%%%%%%%%%%%%%%%%%%%%%%%%%%%%%%%%%%%%%
%%%%%%%%%%%%%%%%%%%%%%%%%%%%%%%%%%%%%%%%%%%%%%%%%%%%%%%%%%%%%%%%%%%%%%%%%%%%%%%%%%%%%%%%%%%%%%%%%%%%%%%%%%%%%%%%%%%%%%%%%%%%%%%%%%%%%%%
%%%%%%%%%%%%%%%%%%%%%%%%%%%%%%%%%%%%%%%%%%%%%%%%%%%%%%%%%%%%%%%%%%%%%%%%%%%%%%%%%%%%%%%%%%%%%%%%%%%%%%%%%%%%%%%%%%%%%%%%%%%%%%%%%%%%%%%

\begin{slide}

\newpage We have the following $ \textbf{natural problems}$ :

\underline{\bf Problem 1}\,: Find the spectral decomposition of
$$L^2(\G_*^J\ba \bHnm)$$ for the Laplacian $\Delta_{n,m}$ on
$\bHnm$ or a commuting set $\BD_*$ of $G^J$-invariant differential
operators on $\bHnm$.

\underline{\bf Problem 2}\,: Decompose the regular representation
$R_{\G_*^J}$ of $G^J$ on $L^2(\G_*^J\ba G^J)$ into irreducibles.

\vskip 1cm \ \ \ The above problems are very important
arithmetically and geometrically. However the above problems are
very {\large\bf difficult} to solve at this moment. One of the
reason is that it is difficult to deal with $\G_*$. Unfortunately
the unitary dual of $Sp(n,\BR)$ is not known yet for $n\geq 3$.

\end{slide}

%%%%%%%%%%%%%%%%%%%%%%%%%%%%%%%%%%%%%%%%%%%%%%%%%%%%%%%%%%%%%%%%%%%%%%%%%%%%%%%%%%%%%%%%%%%%%%%%%%%%%%%%%%%%%%%%%%%%%%%%%%%%%%%%%%%%%%%
%%%%%%%%%%%%%%%%%%%%%%%%%%%%%%%%%%%%%%%%%%%%%%%%%%%%%%%%%%%%%%%%%%%%%%%%%%%%%%%%%%%%%%%%%%%%%%%%%%%%%%%%%%%%%%%%%%%%%%%%%%%%%%%%%%%%%%%
%%%%%%%%%%%%%%%%%%%%%%%%%%%%%%%%%%%%%%%%%%%%%%%%%%%%%%%%%%%%%%%%%%%%%%%%%%%%%%%%%%%%%%%%%%%%%%%%%%%%%%%%%%%%%%%%%%%%%%%%%%%%%%%%%%%%%%%

%\begin{itemize}
%\item $G$ = a Lie group of finite dimension
%\item $H$ = a closed subgroup of $G$
%\vskip 0.3cm \ Then we have a natural representation \\
%$\pi_H$ of $G$ on $L^2(G/H)\,:$
%$$\left( \pi_H(g)f\right)(x):=f(g^{-1}x),$$
%where $g\in G,\ \ f\in L^2(G/H),\ \text{and}\ x\in G/H.$
%\end{itemize}

%\boxed{A\approx B}
%\begin{center}
%\framebox[7cm][c]{ $M_0=M(2,\cO_X)$}
%\end{center}

%%%%%%%%%%%%%%%%%%%%%%%%%%%%%%%%%%%%%%%%%%%%%%%%%%%%%%%%%%%%%%%%%%%%%%%%%%%%%%%%%%%%%%%%%%%%%%%%%%%%%%%%%%%%%%%%%%%%%%%%%%%%%%%%%%%%%%%
%%%%%%%%%%%%%%%%%%%%%%%%%%%%%%%%%%%%%%%%%%%%%%%%%%%%%%%%%%%%%%%%%%%%%%%%%%%%%%%%%%%%%%%%%%%%%%%%%%%%%%%%%%%%%%%%%%%%%%%%%%%%%%%%%%%%%%%
%%%%%%%%%%%%%%%%%%%%%%%%%%%%%%%%%%%%%%%%%%%%%%%%%%%%%%%%%%%%%%%%%%%%%%%%%%%%%%%%%%%%%%%%%%%%%%%%%%%%%%%%%%%%%%%%%%%%%%%%%%%%%%%%%%%%%%%

\begin{slide}

\newpage

For a coordinate $(\O,Z)\in\bHnm$ with $\O=(\o_{\mu\nu})\in {\Bbb
H}_n$ and $Z=(z_{kl})\in \Cmn,$ we put
\begin{align*}
\O=&X+iY,\ \ X=(x_{\mu\nu}),\ \ Y=(y_{\mu\nu})
\  \text{real},\\
Z=&U+iV,\ \ U=(u_{kl}), \ \ V=(v_{kl})\
\text{real},\\
d\O=&(d\o_{\mu\nu}),\ \, d{\overline\O}=(d{\overline \o}_{\mu\nu}),\\
dZ=&(dz_{kl}),\ \ d{\overline Z}=(d{\overline z}_{kl}),
\end{align*}
\begin{align*} \PO=&\left(\, { {1+\delta_{\mu\nu}} \over 2}\, {
{\partial}\over {\partial \o_{\mu\nu}} } \,\right),\ \
\POB=\left(\, { {1+\delta_{\mu\nu}}\over 2} \, { {\partial}\over
{\partial {\overline \o}_{\mu\nu} }  }
\,\right),\\
\PX=&\left(\, { {1+\delta_{\mu\nu}}\over 2}\, { {\partial}\over
{\partial x_{\mu\nu} } } \,\right),\ \ \PY=\left(\, {
{1+\delta_{\mu\nu}}\over 2}\, { {\partial}\over {\partial
y_{\mu\nu} }  } \,\right),
\end{align*}
$$\PZ=\begin{pmatrix} { {\partial}\over{\partial z_{11}} } & \hdots &
{ {\partial}\over{\partial z_{m1}} }\\
\vdots&\ddots&\vdots\\
{ {\partial}\over{\partial z_{1n}} }&\hdots &{ {\partial}\over
{\partial z_{mn}} } \end{pmatrix},$$   $$\PZB=\begin{pmatrix} {
{\partial}\over{\partial {\overline z}_{11} }   }&
\hdots&{ {\partial}\over{\partial {\overline z}_{m1} }  }\\
\vdots&\ddots&\vdots\\
{ {\partial}\over{\partial{\overline z}_{1n} }  }&\hdots & {
{\partial}\over{\partial{\overline z}_{mn} }  }
\end{pmatrix},$$

\end{slide}

%%%%%%%%%%%%%%%%%%%%%%%%%%%%%%%%%%%%%%%%%%%%%%%%%%%%%%%%%%%%%%%%%%%%%%%%%%%%%%%%%%%%%%%%%%%%%%%%%%%%%%%%%%%%%%%%%%%%%%%%%%%%%%%%%%%%%%%
%%%%%%%%%%%%%%%%%%%%%%%%%%%%%%%%%%%%%%%%%%%%%%%%%%%%%%%%%%%%%%%%%%%%%%%%%%%%%%%%%%%%%%%%%%%%%%%%%%%%%%%%%%%%%%%%%%%%%%%%%%%%%%%%%%%%%%%
%%%%%%%%%%%%%%%%%%%%%%%%%%%%%%%%%%%%%%%%%%%%%%%%%%%%%%%%%%%%%%%%%%%%%%%%%%%%%%%%%%%%%%%%%%%%%%%%%%%%%%%%%%%%%%%%%%%%%%%%%%%%%%%%%%%%%%%

\begin{slide}

$$\PU=
\begin{pmatrix} { {\partial}\over{\partial u_{11}} }&\hdots&
{ {\partial}\over{\partial u_{m1}} } \\
\vdots &\ddots &\vdots\\
{ {\partial}\over{\partial u_{1n}} } &\hdots & {
{\partial}\over{\partial u_{mn}} }\end{pmatrix},$$ $$ \PV=
\begin{pmatrix} { {\partial}\over{\partial v_{11}} } &\hdots &
{ {\partial}\over{\partial v_{m1}} }\\
\vdots &\ddots &\hdots\\
{ {\partial}\over{\partial v_{1n}} } &\hdots & {
{\partial}\over{\partial v_{mn}} }\end{pmatrix}.$$

\vskip 3cm

\begin{center}
{\large \bf 1. Invariant metrics on $\bHnm$}
\end{center}

We recall that for a positive real number $A$, the metric
\begin{center}
\framebox[11.5cm][c] {$ds_{n;A}^2=A\cdot \text{tr}\big(Y^{-1}d\O
Y^{-1}d{\overline \O}\big)$}
\end{center}
is a $Sp(n,\BR)$-invariant K{\" a}hler metric on ${\mathbf H}_n$
introduced by C. L. Siegel (cf. [A] or [8], 1943).

\end{slide}

%%%%%%%%%%%%%%%%%%%%%%%%%%%%%%%%%%%%%%%%%%%%%%%%%%%%%%%%%%%%%%%%%%%%%%%%%%%%%%%%%%%%%%%%%%%%%%%%%%%%%%%%%%%%%%%%%%%%%%%%%%%%%%%%%%%%%%%
%%%%%%%%%%%%%%%%%%%%%%%%%%%%%%%%%%%%%%%%%%%%%%%%%%%%%%%%%%%%%%%%%%%%%%%%%%%%%%%%%%%%%%%%%%%%%%%%%%%%%%%%%%%%%%%%%%%%%%%%%%%%%%%%%%%%%%%
%%%%%%%%%%%%%%%%%%%%%%%%%%%%%%%%%%%%%%%%%%%%%%%%%%%%%%%%%%%%%%%%%%%%%%%%%%%%%%%%%%%%%%%%%%%%%%%%%%%%%%%%%%%%%%%%%%%%%%%%%%%%%%%%%%%%%%%

\begin{slide}

\newpage
\noindent {\bf Theorem\ 1 (J.-H. Yang [16], 2005).} For any two
positive real numbers $A$ and $B$, the following metric
\begin{eqnarray*}
& & ds_{n,m;A,B}^2\\
&=&\,A\cdot {\textrm{tr}}\Big( Y^{-1}d\Omega\,Y^{-1}d{\overline \Omega}\Big)
\nonumber
\\ && \ \ + \,B\cdot\bigg\{ {\textrm{tr}}\Big(
Y^{-1}\,{}^tV\,V\,Y^{-1}d\Omega\,Y^{-1} \bo \Big)
\\ &&\quad \quad\quad +\,\,{\textrm{tr}}\Big( Y^{-1}\,{}^t\!(dZ)\,\bz\Big)\\
&&\quad\quad\quad -\,\,{\textrm{tr}}\Big( V\,Y^{-1}d\Om\,Y^{-1}\,{}^t\!(\bz)\Big)\\
&&\quad \quad\quad -\,\,{\textrm{tr}}\Big( V\,Y^{-1}\bo\,
Y^{-1}\,{}^t\!(dZ)\,\Big) \bigg\} \nonumber
\end{eqnarray*}
\noindent is a Riemannian metric on $\bHnm$ which is invariant
under the action (2) of $G^J.$

\end{slide}

%%%%%%%%%%%%%%%%%%%%%%%%%%%%%%%%%%%%%%%%%%%%%%%%%%%%%%%%%%%%%%%%%%%%%%%%%%%%%%%%%%%%%%%%%%%%%%%%%%%%%%%%%%%%%%%%%%%%%%%%%%%%%%%%%%%%%%%
%%%%%%%%%%%%%%%%%%%%%%%%%%%%%%%%%%%%%%%%%%%%%%%%%%%%%%%%%%%%%%%%%%%%%%%%%%%%%%%%%%%%%%%%%%%%%%%%%%%%%%%%%%%%%%%%%%%%%%%%%%%%%%%%%%%%%%%
%%%%%%%%%%%%%%%%%%%%%%%%%%%%%%%%%%%%%%%%%%%%%%%%%%%%%%%%%%%%%%%%%%%%%%%%%%%%%%%%%%%%%%%%%%%%%%%%%%%%%%%%%%%%%%%%%%%%%%%%%%%%%%%%%%%%%%%

\begin{slide}

\vskip 0.2 cm For the case $n=m=A=B=1$, we get
\begin{eqnarray*} & & ds_{1,1;1,1}^2\\ &=&\,{{y\,+\,v^2}\over
{y^3}}\,\big(\,dx^2\,+\,dy^2\,\big)\, +\, {\frac 1y}\,\big(\,du^2\,+\,dv^2\,\big)\\
& &\ \ -\,{{2v}\over {y^2}}\, \big(\,dx\,du\,+\,dy\,dv\,\big).
\end{eqnarray*}

\vskip 0.5cm \noindent {\bf Lemma\ A.} The following differential
form
\begin{eqnarray*}
 dv_{n,m}
={{[dX]\w [dY]\w [dU]\w [dV]}\over {\left(\,\det
Y\,\right)^{n+m+1} } }
\end{eqnarray*}
\noindent is a $G^J$-invariant volume element on $\bHnm$, where
\begin{eqnarray*}
&& [dX]=\w_{\mu\leq\nu}\,dx_{\mu\nu},\quad [dY]=\w_{\mu\leq\nu}
\,dy_{\mu\nu},\\ && [dU]=\w_{k,l}\,du_{kl},\quad
[dV]=\w_{k,l}\,dv_{kl}.
\end{eqnarray*}

$ \textit{Proof.}$ The proof follows from the fact that
\begin{equation*}
(\det Y)^{-(n+1)} [dX]\w [dY]
\end{equation*}
\noindent is a $Sp(n,\BR)$-invariant volume element on $\bH_n$.
(cf. [9])\hfill $\square$

\end{slide}

%%%%%%%%%%%%%%%%%%%%%%%%%%%%%%%%%%%%%%%%%%%%%%%%%%%%%%%%%%%%%%%%%%%%%%%%%%%%%%%%%%%%%%%%%%%%%%%%%%%%%%%%%%%%%%%%%%%%%%%%%%%%%%%%%%%%%%%
%%%%%%%%%%%%%%%%%%%%%%%%%%%%%%%%%%%%%%%%%%%%%%%%%%%%%%%%%%%%%%%%%%%%%%%%%%%%%%%%%%%%%%%%%%%%%%%%%%%%%%%%%%%%%%%%%%%%%%%%%%%%%%%%%%%%%%%
%%%%%%%%%%%%%%%%%%%%%%%%%%%%%%%%%%%%%%%%%%%%%%%%%%%%%%%%%%%%%%%%%%%%%%%%%%%%%%%%%%%%%%%%%%%%%%%%%%%%%%%%%%%%%%%%%%%%%%%%%%%%%%%%%%%%%%%

\begin{slide}

\newpage
\begin{center}
{\large \bf 2. Laplacians on $\bHnm$}
\end{center}

\ \ Hans Maass(cf.\,[3], 1953) proved that for a positive real
number $A$, the differential operator
\begin{center}
\framebox[13cm][c]{$\Delta_n={\frac 4A}\cdot{\textrm{tr}}\left(\,Y
\,{}^{{}^\text{\scriptsize $t$}}\!\!
\left(Y\POB\right)\PO\,\right)$}
\end{center}
is the Laplacian of ${\bH}_n$ for the metric $ds_{n;A}^2.$

[3] H. Maass, {\em Die Differentialgleichungen in der Theorie der
Siegelschen Modulfunktionen}, Math. Ann. {\bf 26} (1953), 44--68.

\end{slide}

%%%%%%%%%%%%%%%%%%%%%%%%%%%%%%%%%%%%%%%%%%%%%%%%%%%%%%%%%%%%%%%%%%%%%%%%%%%%%%%%%%%%%%%%%%%%%%%%%%%%%%%%%%%%%%%%%%%%%%%%%%%%%%%%%%%%%%%
%%%%%%%%%%%%%%%%%%%%%%%%%%%%%%%%%%%%%%%%%%%%%%%%%%%%%%%%%%%%%%%%%%%%%%%%%%%%%%%%%%%%%%%%%%%%%%%%%%%%%%%%%%%%%%%%%%%%%%%%%%%%%%%%%%%%%%%
%%%%%%%%%%%%%%%%%%%%%%%%%%%%%%%%%%%%%%%%%%%%%%%%%%%%%%%%%%%%%%%%%%%%%%%%%%%%%%%%%%%%%%%%%%%%%%%%%%%%%%%%%%%%%%%%%%%%%%%%%%%%%%%%%%%%%%%

\begin{slide}
\newpage
\noindent {\bf Theorem\ 2 (J.-H. Yang [16], 2005).} For any two
positive real numbers $A$ and $B$, the Laplacian
$\Delta_{n,m;A,B}$ of $ds_{n,m;A,B}^2$ is given by
\begin{eqnarray*}
&& \Delta_{n,m;A,B}\\&=& \frac4A \,\bigg\{ {\textrm{tr}}\left(\,Y
\,{}^{{}^\text{\scriptsize $t$}}\!\!
\left(Y\POB\right) \PO\,\right)\\ && \quad\ \ \ +\,{\textrm{tr}}\left(\,VY^{-1}
\,{}^tV   \,{}^{{}^\text{\scriptsize $t$}}\!\!\left(Y\PZB\right)\,\PZ\,\right)\nonumber\\
& &\quad\ \ \
+\,{\textrm{tr}}\left(V
\,{}^{{}^\text{\scriptsize $t$}}\!\!
\left(Y\POB\right) \PZ\,\right)\\ &&\quad \ \ \ +\,{\textrm{tr}}\left(\,{}^tV
\,{}^{{}^\text{\scriptsize $t$}}\!\!
\left(Y\PZB\right) \PO\,\right)\bigg\}\\
& &  +\frac4B\,\,{\textrm{tr}}\left(\, Y\,\PZ
\,{}^{{}^\text{\scriptsize $t$}}\!\!
\left(
\PZB\right)\,\right).\nonumber
\end{eqnarray*}

\end{slide}

%%%%%%%%%%%%%%%%%%%%%%%%%%%%%%%%%%%%%%%%%%%%%%%%%%%%%%%%%%%%%%%%%%%%%%%%%%%%%%%%%%%%%%%%%%%%%%%%%%%%%%%%%%%%%%%%%%%%%%%%%%%%%%%%%%%%%%%
%%%%%%%%%%%%%%%%%%%%%%%%%%%%%%%%%%%%%%%%%%%%%%%%%%%%%%%%%%%%%%%%%%%%%%%%%%%%%%%%%%%%%%%%%%%%%%%%%%%%%%%%%%%%%%%%%%%%%%%%%%%%%%%%%%%%%%%
%%%%%%%%%%%%%%%%%%%%%%%%%%%%%%%%%%%%%%%%%%%%%%%%%%%%%%%%%%%%%%%%%%%%%%%%%%%%%%%%%%%%%%%%%%%%%%%%%%%%%%%%%%%%%%%%%%%%%%%%%%%%%%%%%%%%%%%
%%%%%%%%%%%%%%%%%%%%%%%%%%%%%%%%%%%%%%%%%%%%%%%%%%%%%%%%%%%%%%%%%%%%%%%%%%%%%%%%%%%%%%%%%%%%%%%%%%%%%%%%%%%%%%%%%%%%%%%%%%%%%%%%%%%%%%%

\newcommand\ddx{{{\partial^2}\over{\partial x^2}}}
\newcommand\ddy{{{\partial^2}\over{\partial y^2}}}
\newcommand\ddu{{{\partial^2}\over{\partial u^2}}}
\newcommand\ddv{{{\partial^2}\over{\partial v^2}}}
\newcommand\px{{{\partial}\over{\partial x}}}
\newcommand\py{{{\partial}\over{\partial y}}}
\newcommand\pu{{{\partial}\over{\partial u}}}
\newcommand\pv{{{\partial}\over{\partial v}}}
\newcommand\pxu{{{\partial^2}\over{\partial x\partial u}}}
\newcommand\pyv{{{\partial^2}\over{\partial y\partial v}}}

\begin{slide}

For the case $n=m=A=B=1,$ we get
\begin{eqnarray*}  \Delta_{1,1;1,1}&=&\, y^2\,\left(\,\ddx\,+\,\ddy\,\right)\\ &&\ +\,
(\,y\,+\,v^2\,)\,\left(\,\ddu\,+\,\ddv\,\right)\\ &&\ \,
+\,2\,y\,v\,\left(\,\pxu\,+\,\pyv\,\right).
\end{eqnarray*}

\noindent $ \textbf{Remark\,:}$ $ds^2_{n,m;A,B}$ and
$\Delta_{n,m;A,B}$ are expressed in terms of the $ \textbf{trace
form}$. !!!

\vskip 2cm
\begin{center}
{\large \bf 3. Invariant differential operators on $\bHnm$}
\end{center}

\ \ Let ${\Bbb D}(\bH_n)$ be the algebra of all
$Sp(n,\BR)$-invariant differential operators on $\bH_n$. For
brevity, we set $K=U(n).$ Then $K$ acts on the vector space
\begin{equation*}
T_n=\,\left\{\,\o\in \BC^{(n,n)}\,|\ \o=\,{}^t\o\ \right\}
\end{equation*}

\end{slide}

%%%%%%%%%%%%%%%%%%%%%%%%%%%%%%%%%%%%%%%%%%%%%%%%%%%%%%%%%%%%%%%%%%%%%%%%%%%%%%%%%%%%%%%%%%%%%%%%%%%%%%%%%%%%%%%%%%%%%%%%%%%%%%%%%%%%%%%
%%%%%%%%%%%%%%%%%%%%%%%%%%%%%%%%%%%%%%%%%%%%%%%%%%%%%%%%%%%%%%%%%%%%%%%%%%%%%%%%%%%%%%%%%%%%%%%%%%%%%%%%%%%%%%%%%%%%%%%%%%%%%%%%%%%%%%%
%%%%%%%%%%%%%%%%%%%%%%%%%%%%%%%%%%%%%%%%%%%%%%%%%%%%%%%%%%%%%%%%%%%%%%%%%%%%%%%%%%%%%%%%%%%%%%%%%%%%%%%%%%%%%%%%%%%%%%%%%%%%%%%%%%%%%%%

\begin{slide}

by
\begin{equation}
\framebox[13cm][c]{ $k\cdot \o=\,k\,\o\,{}^tk,\quad h\in K,\ \o\in
T_n .$}
\end{equation}

The action (3) induces naturally the representation $\tau_K$ of $K$ on the
polynomial algebra $ \textrm{Pol}(T_n)$ of $T_n.$ Let
\begin{equation*}
\textrm{Pol}(T_n)^K=\big\{ p\in \textrm{Pol}(T_n)\,\big|\ k\cdot p
=p,\ \forall\, k\in K\,\big\}
\end{equation*}
be the subalgebra of $\textrm{Pol}(T_n)$ consisting of all
$K$-invariant polynomials on $T_n$. Then we get a canonical linear
bijection (not an algebra isomorphism)
\begin{equation}
\mathfrak{S}_n :\textrm{Pol}(T_n)^K\lrt {\Bbb D}(\bH_n).
\end{equation}

\noindent {\bf Theorem\ 3.} $\textrm{Pol}(T_n)^K$ is generated by
algebraically independent polynomials
\begin{equation*}
q_i(\o)=\,{\textrm{tr}}\big( (\o\overline{\o})^i\big),\quad i=1,2,\cdots,n.
\end{equation*}

\end{slide}

%%%%%%%%%%%%%%%%%%%%%%%%%%%%%%%%%%%%%%%%%%%%%%%%%%%%%%%%%%%%%%%%%%%%%%%%%%%%%%%%%%%%%%%%%%%%%%%%%%%%%%%%%%%%%%%%%%%%%%%%%%%%%%%%%%%%%%%
%%%%%%%%%%%%%%%%%%%%%%%%%%%%%%%%%%%%%%%%%%%%%%%%%%%%%%%%%%%%%%%%%%%%%%%%%%%%%%%%%%%%%%%%%%%%%%%%%%%%%%%%%%%%%%%%%%%%%%%%%%%%%%%%%%%%%%%
%%%%%%%%%%%%%%%%%%%%%%%%%%%%%%%%%%%%%%%%%%%%%%%%%%%%%%%%%%%%%%%%%%%%%%%%%%%%%%%%%%%%%%%%%%%%%%%%%%%%%%%%%%%%%%%%%%%%%%%%%%%%%%%%%%%%%%%

\begin{slide}

\noindent $ \textit{Proof.}$ The proof follows from the classical
invariant theory or the work of Harish-Chandra (1923-1983).\hfill
$\square$

\noindent $ \textbf{Remark.}$ Let $D_i=\mathfrak{S}_n(q_i),\ 1\leq
i\leq n$. According to the work of Harish-Chandra,
\begin{equation*}
{\Bbb D}(\bH_n)\cong \BC [D_1,\cdots,D_n]
\end{equation*}
\noindent is a polynomial ring of degree $n$, where $n$ is the
split real rank of $Sp(n,\BR).$

\noindent $ \textbf{Remark.}$ $\mathfrak{S}_n(q_1)=\Delta_{n;1}$
is the Laplacian of $ds_{n;1}^2$ on $\bH_n.$ So far
$\mathfrak{S}_n(q_i)\,(i=2,\cdots,n)$ were not written explicitly.

\noindent $ \textbf{Remark.}$
Maass [3] found explicit algebraically independent generators $H_1,H_2,\cdots,H_n$ of ${\mathbb D}(\BH_n)$.
We will describe $H_1,H_2,\cdots,H_n$ explicitly. For $M=\begin{pmatrix} A&B\\
C&D\end{pmatrix} \in Sp(n,\BR)$ and $\Omega=X+iY\in \BH_n$ with real $X,Y$, we set
\begin{equation*}
\Omega_*=\,M\!\cdot\!\Omega=\,X_*+\,iY_*\quad \textrm{with}\ X_*,Y_*\ \textrm{real}.
\end{equation*}

\end{slide}

%%%%%%%%%%%%%%%%%%%%%%%%%%%%%%%%%%%%%%%%%%%%%%%%%%%%%%%%%%%%%%%%%%%%%%%%%%%%%%%%%%%%%%%%%%%%%%%%%%%%%%%%%%%%%%%%%%%%%%%%%%%%%%%%%%%%%%%
%%%%%%%%%%%%%%%%%%%%%%%%%%%%%%%%%%%%%%%%%%%%%%%%%%%%%%%%%%%%%%%%%%%%%%%%%%%%%%%%%%%%%%%%%%%%%%%%%%%%%%%%%%%%%%%%%%%%%%%%%%%%%%%%%%%%%%%
%%%%%%%%%%%%%%%%%%%%%%%%%%%%%%%%%%%%%%%%%%%%%%%%%%%%%%%%%%%%%%%%%%%%%%%%%%%%%%%%%%%%%%%%%%%%%%%%%%%%%%%%%%%%%%%%%%%%%%%%%%%%%%%%%%%%%%%

\begin{slide}

We set
\begin{eqnarray*}
K&=&\,\big( \Omega-{\overline\Omega}\,\big)\PO=\,2\,i\,Y \PO,\\
\Lambda&=&\,\big( \Omega-{\overline\Omega}\,\big)\POB=\,2\,i\,Y \POB,\\
K_*&=& \,\big( \Omega_*-{\overline\Omega}_*\,\big)\POS=\,2\,i\,Y_* \POS,\\
\Lambda_*&=&\,\big( \Omega_*-{\overline\Omega}_*\,\big)\POBS=\,2\,i\,Y_* \POBS.
\end{eqnarray*}
Then it is easily seen that
\begin{equation}
K_*=\,{}^t(C{\overline\Om}+D)^{-1}\,{}^t\!\left\{ (C\Omega+D)\,{}^t\!K \right\},
\end{equation}

\begin{equation}
\Lambda_*=\,{}^t(C{\Om}+D)^{-1}\,{}^t\!\left\{ (C{\overline\Omega}+D)\,{}^t\!\Lambda \right\}
\end{equation}
and
\begin{equation}
{}^t\!\left\{ (C{\overline\Omega}+D)\,{}^t\!\Lambda \right\}=\,\Lambda\,{}^t(C{\overline\Omega}+D)
-{{n+1}\over 2} \,\big( \Omega-{\overline\Omega}\,\big)\,{}^t\!C.
\end{equation}

Using Formulas (5),\,(6) and (7), we can show that

\end{slide}

%%%%%%%%%%%%%%%%%%%%%%%%%%%%%%%%%%%%%%%%%%%%%%%%%%%%%%%%%%%%%%%%%%%%%%%%%%%%%%%%%%%%%%%%%%%%%%%%%%%%%%%%%%%%%%%%%%%%%%%%%%%%%%%%%%%%%%%
%%%%%%%%%%%%%%%%%%%%%%%%%%%%%%%%%%%%%%%%%%%%%%%%%%%%%%%%%%%%%%%%%%%%%%%%%%%%%%%%%%%%%%%%%%%%%%%%%%%%%%%%%%%%%%%%%%%%%%%%%%%%%%%%%%%%%%%
%%%%%%%%%%%%%%%%%%%%%%%%%%%%%%%%%%%%%%%%%%%%%%%%%%%%%%%%%%%%%%%%%%%%%%%%%%%%%%%%%%%%%%%%%%%%%%%%%%%%%%%%%%%%%%%%%%%%%%%%%%%%%%%%%%%%%%%

\begin{slide}

\begin{eqnarray*}
& & \Lambda_*K_* \,+\,{{n+1}\over 2}K_*\\
&=&\,{}^t(C{\Om}+D)^{-1}\,{}^{{}^{{}^{{}^\text{\scriptsize $t$}}}}\!\!\!
\left\{ (C{\Omega}+D)\,{}^{{}^{{}^{{}^\text{\scriptsize $t$}}}}\!\!\!
\left( \Lambda K \,+\, {{n+1}\over 2}K \right)\right\}.
\end{eqnarray*}

\noindent Therefore we get
\begin{equation*}
\textrm{tr}\!  \left( \Lambda_*K_* \,+\,{{n+1}\over 2}K_* \right) =\,
 \textrm{tr}\!   \left( \Lambda K \,+\, {{n+1}\over 2}K \right).
\end{equation*}
We set
\begin{equation*}
A^{(1)}=\,\Lambda K \,+\, {{n+1}\over 2}K .
\end{equation*}

We define $A^{(j)}\,(j=2,3,\cdots,n)$ recursively by
\begin{eqnarray}
A^{(j)}&=&\, A^{(1)}A^{(j-1)}- {{n+1}\over 2}\,\Lambda\, A^{(j-1)}\nonumber\\
& &\, \,+\,{\frac 12}\,\Lambda\,\sigma\!\left(
A^{(j-1)} \right)\\
& & \, \,+\, {\frac 12}\,\big( \Omega-{\overline\Omega}\,\big) \,
{}^{{}^{{}^\text{\scriptsize $t$}}}\!\!\!\left\{
\big( \Omega-{\overline\Omega}\,\big)^{-1}
\,{}^t\!\left( \,{}^t\!\Lambda\,{}^t\!A^{(j-1)}\right)\right\}.\nonumber
\end{eqnarray}

\end{slide}

%%%%%%%%%%%%%%%%%%%%%%%%%%%%%%%%%%%%%%%%%%%%%%%%%%%%%%%%%%%%%%%%%%%%%%%%%%%%%%%%%%%%%%%%%%%%%%%%%%%%%%%%%%%%%%%%%%%%%%%%%%%%%%%%%%%%%%%
%%%%%%%%%%%%%%%%%%%%%%%%%%%%%%%%%%%%%%%%%%%%%%%%%%%%%%%%%%%%%%%%%%%%%%%%%%%%%%%%%%%%%%%%%%%%%%%%%%%%%%%%%%%%%%%%%%%%%%%%%%%%%%%%%%%%%%%
%%%%%%%%%%%%%%%%%%%%%%%%%%%%%%%%%%%%%%%%%%%%%%%%%%%%%%%%%%%%%%%%%%%%%%%%%%%%%%%%%%%%%%%%%%%%%%%%%%%%%%%%%%%%%%%%%%%%%%%%%%%%%%%%%%%%%%%

\begin{slide}

\noindent We set
\begin{equation}
H_j=\,\textrm{tr}\!  \left( A^{(j)} \right),\quad j=1,2,\cdots,n.
\end{equation}
As mentioned before, Maass proved that $H_1,H_2,\\
\cdots,H_n$ are algebraically independent generators
$H_1,H_2,\cdots,H_n$ of ${\mathbb D}(\BH_n)$.

**************************************
**************************************

\ \ \ Let $T_{n,m}=T_n \times \Cmn.$ Then $K$ acts on $T_{n,m}$ by
%\begin{center}
%\framebox[7cm][c]{ $M_0=M(2,\cO_X)$}
%\end{center}

\begin{equation}
\framebox[11cm][c]{
$h\cdot(\o,z)=\,\big(\,h\,\o\,{}^th,\,z\,{}^th\,\big),$}
\end{equation}
\noindent where $h\in K,\ \o\in T_n,\ z\in \Cmn.$ Then this action
induces naturally the action $\rho$ of $K$ on the polynomial
algebra
$$\text{Pol}_{m,n}=\,\text{Pol}\,(T_{n,m}).$$

\end{slide}

%%%%%%%%%%%%%%%%%%%%%%%%%%%%%%%%%%%%%%%%%%%%%%%%%%%%%%%%%%%%%%%%%%%%%%%%%%%%%%%%%%%%%%%%%%%%%%%%%%%%%%%%%%%%%%%%%%%%%%%%%%%%%%%%%%%%%%%
%%%%%%%%%%%%%%%%%%%%%%%%%%%%%%%%%%%%%%%%%%%%%%%%%%%%%%%%%%%%%%%%%%%%%%%%%%%%%%%%%%%%%%%%%%%%%%%%%%%%%%%%%%%%%%%%%%%%%%%%%%%%%%%%%%%%%%%
%%%%%%%%%%%%%%%%%%%%%%%%%%%%%%%%%%%%%%%%%%%%%%%%%%%%%%%%%%%%%%%%%%%%%%%%%%%%%%%%%%%%%%%%%%%%%%%%%%%%%%%%%%%%%%%%%%%%%%%%%%%%%%%%%%%%%%%

\begin{slide}

We denote by
$\text{Pol}_{m,n}^K$ the subalgebra of $\text{Pol}_{m,n}$
consisting of all $K$-invariants of the action $\rho$ of $K.$ We
also denote by $$\BD(\bHnm)$$ the algebra of all differential
operators on $\bHnm$ which is invariant under the action (2) of
the Jacobi group $G^J$. Then we can show that there exists a
natural linear bijection
$$\mathfrak{S}_{n,m}:\,\text{Pol}^K_{m,n}\lrt \BD(\bHnm)$$
of $\text{Pol}^K_{m,n}$ onto $\BD(\bHnm).$

\vskip 0.2cm \ \ \ The map $\mathfrak{S}_{n,m}$ is described
explicitly as follows. \par We put $N_{\star}=n(n+1)+2mn$. Let
$\big\{ \eta_{\alpha}\,|\ 1\leq \alpha \leq N_{\star}\, \big\}$ be
a basis of $T_{n,m}$. If $P\in \mathrm{Pol}_{m,n}^K$, then
\begin{align*}
& & \Big(\mathfrak{S}_{n,m} (P)f\Big)(gK) \hskip 11cm\\
&=&\left[ P\left( {{\partial}\over {\partial
t_{\alpha}}}\right)f\left(g\,\text{exp}\,
\left(\sum_{\alpha=1}^{N_{\star}} t_{\alpha}\eta_{\alpha}\right)
K\right)\right]_{(t_{\alpha})=0},\hskip 3cm
\end{align*}

\end{slide}

%%%%%%%%%%%%%%%%%%%%%%%%%%%%%%%%%%%%%%%%%%%%%%%%%%%%%%%%%%%%%%%%%%%%%%%%%%%%%%%%%%%%%%%%%%%%%%%%%%%%%%%%%%%%%%%%%%%%%%%%%%%%%%%%%%%%%%%
%%%%%%%%%%%%%%%%%%%%%%%%%%%%%%%%%%%%%%%%%%%%%%%%%%%%%%%%%%%%%%%%%%%%%%%%%%%%%%%%%%%%%%%%%%%%%%%%%%%%%%%%%%%%%%%%%%%%%%%%%%%%%%%%%%%%%%%
%%%%%%%%%%%%%%%%%%%%%%%%%%%%%%%%%%%%%%%%%%%%%%%%%%%%%%%%%%%%%%%%%%%%%%%%%%%%%%%%%%%%%%%%%%%%%%%%%%%%%%%%%%%%%%%%%%%%%%%%%%%%%%%%%%%%%%%

\begin{slide}

\noindent where $f\in C^{\infty}({\mathbb H}_{n,m})$. In general,
it is hard to express $\mathfrak{S}_{n,m}(P)$ explicitly for a
polynomial $P\in \textrm{Pol}_{m,n}^K$.

**************************************
**************************************

\vskip 0.2cm \vskip 0.2cm  We present the following {\bf basic} $K$-invariant
polynomials in $\text{Pol}_{m,n}^K$.
\begin{eqnarray*}
&& p_j(\om,z)=\,\text{tr}((\om{\overline
\om})^j),\quad 1\leq j\leq n,\\
&&  \psi_k^{(1)}(\om,z)=\,(z\,^t{\overline
z})_{kk},\quad 1\leq k\leq m,  \\
&&  \psi_{kp}^{(2)}(\om,z)=
\,\text{Re}\,(z\,^t{\overline z})_{kp}, \quad 1\leq k<p\leq m,\\
&& \psi_{kp}^{(3)}(\om,z)
=\,\text{Im}\,(z\,^t{\overline z})_{kp},\quad 1\leq k<p\leq m,\\
&& f_{kp}^{(1)}(\om,z)= \,\text{Re}\,(z{\overline
\om}\,^tz)_{kp},\quad 1\leq k\leq p\leq m,\\
&& f_{kp}^{(2)}(\om,z)= \,\text{Im}\,(z{\overline
\om}\,^tz\,)_{kp},\quad 1\leq k\leq p\leq m,
\end{eqnarray*}
\noindent where $\om\in T_n$ and $z\in \BC^{(m,n)}$.

\end{slide}

%%%%%%%%%%%%%%%%%%%%%%%%%%%%%%%%%%%%%%%%%%%%%%%%%%%%%%%%%%%%%%%%%%%%%%%%%%%%%%%%%%%%%%%%%%%%%%%%%%%%%%%%%%%%%%%%%%%%%%%%%%%%%%%%%%%%%%%
%%%%%%%%%%%%%%%%%%%%%%%%%%%%%%%%%%%%%%%%%%%%%%%%%%%%%%%%%%%%%%%%%%%%%%%%%%%%%%%%%%%%%%%%%%%%%%%%%%%%%%%%%%%%%%%%%%%%%%%%%%%%%%%%%%%%%%%
%%%%%%%%%%%%%%%%%%%%%%%%%%%%%%%%%%%%%%%%%%%%%%%%%%%%%%%%%%%%%%%%%%%%%%%%%%%%%%%%%%%%%%%%%%%%%%%%%%%%%%%%%%%%%%%%%%%%%%%%%%%%%%%%%%%%%%%

\begin{slide}

\vskip 0.3cm For an $m\times m$ matrix $S$, we define the
following invariant polynomials in $\text{Pol}_{m,n}^K$.
\begin{eqnarray*}
&&
m_{j;S}^{(1)}(\om,z)=\,\textrm{Re}\,\Big(\text{tr}\big(\om{\overline
\om}+ \,^tzS{\overline
z}\,\big)^j\,\Big),\\
&&
m_{j;S}^{(2)}(\om,z)=\,\textrm{Im}\,\Big(\text{tr}\big(\om{\overline
\om}+ \,^tzS{\overline
z}\,\big)^j\,\Big),\\
&&  q_{k;S}^{(1)}(\om,z)=\,\textrm{Re}\,\Big( \textrm{tr}\big( (
\,^tz\,S\,{\overline
z})^k\big) \Big),  \\
&&  q_{k;S}^{(2)}(\om,z)=\,\textrm{Im}\,\Big( \textrm{tr}\big( (
\,^tz\,S\,{\overline
z})^k\big) \Big),  \\
&& \theta_{i,k,j;S}^{(1)}(\om,z) \\
&=&\,\textrm{Re}\,\Big( \textrm{tr}\big( (\om {\overline
\om})^i\,(\,^tz\,S\,{\overline z})^k\,(\om {\overline
\om}+\,^tz\,S\,{\overline
z}\,)^j\,\big)\Big),\\
&& \theta_{i,k,j;S}^{(2)}(\om,z)\\
&=&\,\textrm{Im}\,\Big( \textrm{tr}\big( (\om {\overline
\om})^i\,(\,^tz\,S\,{\overline z})^k\,(\om {\overline
\om}+\,^tz\,S\,{\overline z}\,)^j\,\big)\Big),
\end{eqnarray*}
\noindent where $1\leq i,j\leq n$ and $1\leq k\leq m$.

\end{slide}

%%%%%%%%%%%%%%%%%%%%%%%%%%%%%%%%%%%%%%%%%%%%%%%%%%%%%%%%%%%%%%%%%%%%%%%%%%%%%%%%%%%%%%%%%%%%%%%%%%%%%%%%%%%%%%%%%%%%%%%%%%%%%%%%%%%%%%%
%%%%%%%%%%%%%%%%%%%%%%%%%%%%%%%%%%%%%%%%%%%%%%%%%%%%%%%%%%%%%%%%%%%%%%%%%%%%%%%%%%%%%%%%%%%%%%%%%%%%%%%%%%%%%%%%%%%%%%%%%%%%%%%%%%%%%%%
%%%%%%%%%%%%%%%%%%%%%%%%%%%%%%%%%%%%%%%%%%%%%%%%%%%%%%%%%%%%%%%%%%%%%%%%%%%%%%%%%%%%%%%%%%%%%%%%%%%%%%%%%%%%%%%%%%%%%%%%%%%%%%%%%%%%%%%

\begin{slide}

\newpage
We define the following $K$-invariant polynomials in
$\text{Pol}_{m,n}^K$.
\begin{eqnarray*}
&&  r_{jk}^{(1)}(\om,z)= \,\textrm{Re}\,\Big( \textrm{tr}\big(
(\om {\overline \om})^j\,(\,^tz{\overline z})^k\,\big)\Big),
\\
&&  r_{jk}^{(2)}(\om,z)= \,\textrm{Im}\,\Big( \textrm{tr}\big(
(\om {\overline \om})^j\,(\,^tz{\overline z})^k\,\big)\Big),
\end{eqnarray*}
\noindent where $1\leq j\leq n$ and $1\leq k\leq m$.

\vskip 0.3cm There may be possible other new invariants. We think
that at this moment it may be complicated and difficult to find
the generators of $\text{Pol}_{m,n}^K$.

 \vskip 0.355cm We propose the following problems.

\vskip 0.2cm \noindent $ \textbf{Problem A.}$ Find the generators
of $\text{Pol}_{m,n}^K$.

\vskip 0.32cm \noindent $ \textbf{Problem B.}$ Find an easy way to
express the images of the above invariant polynomials under the
map $\mathfrak{S}_{n,m}$ explicitly.

\end{slide}

%%%%%%%%%%%%%%%%%%%%%%%%%%%%%%%%%%%%%%%%%%%%%%%%%%%%%%%%%%%%%%%%%%%%%%%%%%%%%%%%%%%%%%%%%%%%%%%%%%%%%%%%%%%%%%%%%%%%%%%%%%%%%%%%%%%%%%%
%%%%%%%%%%%%%%%%%%%%%%%%%%%%%%%%%%%%%%%%%%%%%%%%%%%%%%%%%%%%%%%%%%%%%%%%%%%%%%%%%%%%%%%%%%%%%%%%%%%%%%%%%%%%%%%%%%%%%%%%%%%%%%%%%%%%%%%
%%%%%%%%%%%%%%%%%%%%%%%%%%%%%%%%%%%%%%%%%%%%%%%%%%%%%%%%%%%%%%%%%%%%%%%%%%%%%%%%%%%%%%%%%%%%%%%%%%%%%%%%%%%%%%%%%%%%%%%%%%%%%%%%%%%%%%%

%\newpage \noindent {\bf Conjecture.} The algebra
%$\BD(\bHnm)$ is generated by the images under the mapping
%$\mathfrak{S}_{n,m}$ of the following invariants\par\ \ \ \ \ \ \
%$\text{(I1)} \ p_j(\o,z)=\,\text{tr}((\o{\overline \o})^j),$\par\
%\ \ \ \ \ \ $\text{(I2)}\ \psi_k^{(1)}(\o,z)=\,(z\,{\overline
%z}^T\,)_{kk},$\par\ \ \ \ \ \ \ $ \text{(I3)}\
%\psi_{kp}^{(2)}(\o,z)= \,\text{Re}\,(z\,{\overline z}^T\,)_{kp},\
%\ k< p, $\par \ \ \ \ \ \ \ $\text{(I4)}\ \psi_{kp}^{(3)}(\o,z)
%=\,\text{Im}\,(z\,{\overline z}^T\,)_{kp},\ \ k< p, $\par \ \ \ \
%\ \ \ $\text{(I5)}\ \varphi_{kp}^{(1)}(\o,z)=
%\,\text{Re}\,(z{\overline \o}\,z^T)_{kp}, $\par \ \ \ \ \ \ \
%$\text{(I6)}\ \varphi_{kp}^{(2)}(z,w)= \,\text{Im}\,(z{\overline
%\o}\,z^T\,)_{kp}, $\par where $1\leq j\leq n,\ 1\leq k\leq m$ and
%$1\leq k\leq p\leq m$. Moreover, $\BD(\bHnm)$ is not commutative.

\def\ddx{{{\partial^2}\over{\partial x^2}}}
\def\ddy{{{\partial^2}\over{\partial y^2}}}
\def\ddu{{{\partial^2}\over{\partial u^2}}}
\def\ddv{{{\partial^2}\over{\partial v^2}}}
\def\px{{{\partial}\over{\partial x}}}
\def\py{{{\partial}\over{\partial y}}}
\def\pu{{{\partial}\over{\partial u}}}
\def\pv{{{\partial}\over{\partial v}}}
\def\pxu{{{\partial^2}\over{\partial x\partial u}}}
\def\pyv{{{\partial^2}\over{\partial y\partial v}}}
\def\DSPR{{\Bbb D}(\SPR)}
\def\dx{{{\partial}\over{\partial x}}}
\def\dy{{{\partial}\over{\partial y}}}
\def\du{{{\partial}\over{\partial u}}}
\def\dv{{{\partial}\over{\partial v}}}

%%%%%%%%%%%%%%%%%%%%%%%%%%%%%%%%%%%%%%%%%%%%%%%%%%%%%%%%%%%%%%%%%%%%%%%%%%%%%%%%%%%%%%%%%%%%%%%%%%%%%%%%%%%%%%%%%%%%%%%%%%%%%%%%%%%%%%%
%%%%%%%%%%%%%%%%%%%%%%%%%%%%%%%%%%%%%%%%%%%%%%%%%%%%%%%%%%%%%%%%%%%%%%%%%%%%%%%%%%%%%%%%%%%%%%%%%%%%%%%%%%%%%%%%%%%%%%%%%%%%%%%%%%%%%%%
%%%%%%%%%%%%%%%%%%%%%%%%%%%%%%%%%%%%%%%%%%%%%%%%%%%%%%%%%%%%%%%%%%%%%%%%%%%%%%%%%%%%%%%%%%%%%%%%%%%%%%%%%%%%%%%%%%%%%%%%%%%%%%%%%%%%%%%

\begin{slide}

\newpage
{\bf Theorem\ 4.} The algebra $\BD({\mathbf H}_1\times \BC)$ is
generated by the following differential operators \begin{align*}
D=&y^2\,\left(\,{{\partial^2}\over {\partial x^2}}+
{{\partial^2}\over {\partial y^2}}\,\right)
+v^2\,\left(\,\ddu\,+\,\ddv\,\right) \\
&\ \ +2\,y\,v\,\left(\,\pxu\,+\,\pyv\,\right),
\end{align*}
$$\Psi=y\left(\,{{\partial^2}\over {\partial u^2}}+
{{\partial^2}\over {\partial v^2}}\,\right),\hskip 7cm$$
\begin{align*}D_1=&\,2y^2\,{{\partial^3}\over {\partial x\partial u
\partial v}}-y^2\,{{\partial}\over{\partial y}}
\left(\,{{\partial^2}\over{\partial u^2}}-
{{\partial^2}\over{\partial v^2}}\,\right)\hskip 1cm\\ &\ \ \
+\left(\, v\,{{\partial}\over{\partial v}}\,+\,1\,\right)\Psi
\end{align*} and
\begin{align*} D_2=&\,y^2\,{{\partial}\over{\partial x}}\left(\,
{{\partial^2}\over{\partial v^2}}\,-\,{{\partial^2}\over {\partial
u^2}}\,\right)\,-\,2\,y^2\,{{\partial^3}\over{\partial y\partial u
\partial v}}\\ &\ \ \ \ -\,v\,{{\partial}\over{\partial u}}\Psi,
\end{align*} where $\tau=x+iy$ and $z=u+iv$ with real

\end{slide}

%%%%%%%%%%%%%%%%%%%%%%%%%%%%%%%%%%%%%%%%%%%%%%%%%%%%%%%%%%%%%%%%%%%%%%%%%%%%%%%%%%%%%%%%%%%%%%%%%%%%%%%%%%%%%%%%%%%%%%%%%%%%%%%%%%%%%%%
%%%%%%%%%%%%%%%%%%%%%%%%%%%%%%%%%%%%%%%%%%%%%%%%%%%%%%%%%%%%%%%%%%%%%%%%%%%%%%%%%%%%%%%%%%%%%%%%%%%%%%%%%%%%%%%%%%%%%%%%%%%%%%%%%%%%%%%
%%%%%%%%%%%%%%%%%%%%%%%%%%%%%%%%%%%%%%%%%%%%%%%%%%%%%%%%%%%%%%%%%%%%%%%%%%%%%%%%%%%%%%%%%%%%%%%%%%%%%%%%%%%%%%%%%%%%%%%%%%%%%%%%%%%%%%%

\begin{slide}

variables
$x,y,u,v.$ Moreover, we have \begin{align*} D\Psi-&\Psi D\,=\,
2\,y^2\,\dy\left(\,\ddu\,-\,\ddv\,\right)\\ & -
4\,y^2\,{{\partial^3}\over{\partial x\partial u\partial
v}}-2\,\left(\,v\,\dv\Psi+\Psi\,\right).
\end{align*}\par
%In particular, the algebra $\BD({\mathbf H}_1\times \BC)$ is not
%commutative.
\def\bz{d{\overline Z}}
\def\bo{d{\overline \O}}
\def\BZ{\Bbb Z}
\def\s{\text{tr}}

\noindent $ \textbf{Remark.}$ We observe that $\Delta_{n,m;A,B}\in
\BD(\bHnm).$ We can show that
$$D={\textrm{tr}}\left(\, Y\,\PZ
{}^{{}^{{}^\text{\scriptsize $t$}}}\!\!\!
\left(
\PZB\right)\,\right)$$ \noindent is an element of $\BD(\bHnm).$
Therefore
\begin{equation*}
\Delta_{n,m;A,B}- {4 \over B}\,D \in \BD(\bHnm).
\end{equation*}
The following differential
operator ${\mathbb K}$ on ${\mathbb H}_{n,m}$ of degree $2n$ defined by
\begin{equation*}
{\mathbb K}=\,\det(Y)\,\det\left( \PZ {}^{{}^{{}^\text{\scriptsize $t$}}}\!\!\!\left(
\PZB\right)\right)
\end{equation*}
is invariant under the action (2) of $G^J$.

\end{slide}

%%%%%%%%%%%%%%%%%%%%%%%%%%%%%%%%%%%%%%%%%%%%%%%%%%%%%%%%%%%%%%%%%%%%%%%%%%%%%%%%%%%%%%%%%%%%%%%%%%%%%%%%%%%%%%%%%%%%%%%%%%%%%%%%%%%%%%%
%%%%%%%%%%%%%%%%%%%%%%%%%%%%%%%%%%%%%%%%%%%%%%%%%%%%%%%%%%%%%%%%%%%%%%%%%%%%%%%%%%%%%%%%%%%%%%%%%%%%%%%%%%%%%%%%%%%%%%%%%%%%%%%%%%%%%%%
%%%%%%%%%%%%%%%%%%%%%%%%%%%%%%%%%%%%%%%%%%%%%%%%%%%%%%%%%%%%%%%%%%%%%%%%%%%%%%%%%%%%%%%%%%%%%%%%%%%%%%%%%%%%%%%%%%%%%%%%%%%%%%%%%%%%%%%

\begin{slide}

\newpage
The following
matrix-valued differential operator ${\mathbb T}$ on
$\BH_{n,m}$ defined by
\begin{equation*}
{\mathbb T}=\,
{}^{{}^{{}^{{}^\text{\scriptsize $t$}}}}\!\!\!
\left( \PZB\right) Y \PZ
\end{equation*}

\noindent is invariant under the action (2) of $G^J$.
Therefore each $(k,l)$-entry ${\mathbb T}_{kl}$ of ${\mathbb
T}$ given by
\begin{equation*}
{\mathbb T}_{kl}=\sum_{i,j=1}^n
\,y_{ij}\,{{\partial^2\ \ \ \ }\over{\partial {\overline z}_{ki}\partial
z_{lj}} },\quad 1\leq k,l\leq m
\end{equation*}

\noindent is an element of $\BD\big(\BH_{n,m}\big)$.

\vskip 0.1cm
Indeed
it is very complicated and difficult at this moment to express the
generators of the algebra of all $G^J$-invariant differential
operators on $\Dnm$ explicitly. In particular, it is extremely difficult to
find explicit $G^J$-invariant differential operators on
${\mathbb H}_{n,m}$ of {\it odd} degree.
We propose an open problem to find
other explicit $G^J$-invariant differential operators on
${\mathbb H}_{n,m}$.

\end{slide}

%%%%%%%%%%%%%%%%%%%%%%%%%%%%%%%%%%%%%%%%%%%%%%%%%%%%%%%%%%%%%%%%%%%%%%%%%%%%%%%%%%%%%%%%%%%%%%%%%%%%%%%%%%%%%%%%%%%%%%%%%%%%%%%%%%%%%%%
%%%%%%%%%%%%%%%%%%%%%%%%%%%%%%%%%%%%%%%%%%%%%%%%%%%%%%%%%%%%%%%%%%%%%%%%%%%%%%%%%%%%%%%%%%%%%%%%%%%%%%%%%%%%%%%%%%%%%%%%%%%%%%%%%%%%%%%
%%%%%%%%%%%%%%%%%%%%%%%%%%%%%%%%%%%%%%%%%%%%%%%%%%%%%%%%%%%%%%%%%%%%%%%%%%%%%%%%%%%%%%%%%%%%%%%%%%%%%%%%%%%%%%%%%%%%%%%%%%%%%%%%%%%%%%%

\begin{slide}

\newpage
\begin{center}
{\large \bf 4. Partial Cayley transform}
\end{center}

\ \ Let $$\bDn=\left\{W\in\BC^{(n,n)}\,|\ W=\,{}^tW,\ I_n-W{\overline
W}>0\right\}$$ be the generalized unit disk of degree $n$. We let
\begin{equation*}
\bDnm=\bDn\times \Cmn
\end{equation*}
\noindent be the Siegel-Jacobi disk. \par \ \ We define the $
\textbf{partial Cayley transform}$
\begin{equation*}
\Phi_*:\bDnm\lrt \bHnm
\end{equation*}
by \noindent
\begin{equation}
\Phi_*(W,\eta)=
\end{equation}
\begin{equation*}
\framebox[15cm][c]{$\left(
i(I_n+W)(I_n-W)^{-1},\,2\,i\,\eta\,(I_n-W)^{-1}\right),$}
\end{equation*}
\noindent where $W\in \bDn$ and $\eta\in \Cmn.$ It is easy to see
that $\Phi_*$ is a biholomorphic mapping.

\end{slide}

%%%%%%%%%%%%%%%%%%%%%%%%%%%%%%%%%%%%%%%%%%%%%%%%%%%%%%%%%%%%%%%%%%%%%%%%%%%%%%%%%%%%%%%%%%%%%%%%%%%%%%%%%%%%%%%%%%%%%%%%%%%%%%%%%%%%%%%
%%%%%%%%%%%%%%%%%%%%%%%%%%%%%%%%%%%%%%%%%%%%%%%%%%%%%%%%%%%%%%%%%%%%%%%%%%%%%%%%%%%%%%%%%%%%%%%%%%%%%%%%%%%%%%%%%%%%%%%%%%%%%%%%%%%%%%%
%%%%%%%%%%%%%%%%%%%%%%%%%%%%%%%%%%%%%%%%%%%%%%%%%%%%%%%%%%%%%%%%%%%%%%%%%%%%%%%%%%%%%%%%%%%%%%%%%%%%%%%%%%%%%%%%%%%%%%%%%%%%%%%%%%%%%%%

\begin{slide}

\vskip 0.1cm We set
\begin{equation*}
T_*={1\over {\sqrt 2}}\,
\begin{pmatrix} I_{m+n} & I_{m+n}\\ iI_{m+n} & -iI_{m+n}
\end{pmatrix}.
\end{equation*}
We now consider the group $G_*^J$ defined by
\begin{equation*}
G_*^J=T_*^{-1}G^JT_*.
\end{equation*}
Then $G_*^J$ acts on $\bDnm$ transitively by
\begin{equation}
 \left(\begin{pmatrix} P & Q\\
{\overline Q} & {\overline P}
\end{pmatrix},\left( \la, \mu,\kappa\right)\right)\cdot
(W,\eta)=
\end{equation}
\begin{equation*}
\big((PW+Q)(\OQ W+\OP)^{-1},(\eta+\la W+\mu)(\OQ
W+\OP)^{-1}\big).\notag\end{equation*}

\vskip 0.3cm\noindent {\bf Theorem 5 (J.-H. Yang [17], 2005).} The
action (2) of $G^J$ on $\bHnm$ is compatible with the action (12)
of $G_*^J$ on $\bDnm$ through the partial Cayley transform
$\Phi_*$. More precisely, if $g_0\in G^J$ and $(W,\eta)\in\bDnm$,
\begin{equation*}
g_0\cdot\Phi_*(W,\eta)=\Phi_*\big(g_*\cdot (W,\eta)\big),
\end{equation*}
where $g_*= T_*^{-1}g_0 T_*$.

\end{slide}

%%%%%%%%%%%%%%%%%%%%%%%%%%%%%%%%%%%%%%%%%%%%%%%%%%%%%%%%%%%%%%%%%%%%%%%%%%%%%%%%%%%%%%%%%%%%%%%%%%%%%%%%%%%%%%%%%%%%%%%%%%%%%%%%%%%%%%%
%%%%%%%%%%%%%%%%%%%%%%%%%%%%%%%%%%%%%%%%%%%%%%%%%%%%%%%%%%%%%%%%%%%%%%%%%%%%%%%%%%%%%%%%%%%%%%%%%%%%%%%%%%%%%%%%%%%%%%%%%%%%%%%%%%%%%%%
%%%%%%%%%%%%%%%%%%%%%%%%%%%%%%%%%%%%%%%%%%%%%%%%%%%%%%%%%%%%%%%%%%%%%%%%%%%%%%%%%%%%%%%%%%%%%%%%%%%%%%%%%%%%%%%%%%%%%%%%%%%%%%%%%%%%%%%

\newcommand\bw{d{\overline W}}
\newcommand\be{d{\overline \eta}}

\begin{slide}

\newpage
\begin{center}
{\large \bf 5. Invariant Differential Operators on $\bDnm$}
\end{center}
\ \ For a coordinate $(W,\eta)\in \bDnm$ with $W=(w_{\mu\nu})\in
{\mathbf D}_n$ and $\eta=(\eta_{kl})\in \Cmn,$ we put
\begin{eqnarray*}
dW\,&=&\,(dw_{\mu\nu}),\quad\ \ d{\overline W}\,=\,(d{\overline w}_{\mu\nu}),\\
d\eta\,&=&\,(d\eta_{kl}),\quad\ \
d{\overline\eta}\,=\,(d{\overline\eta}_{kl}),
\end{eqnarray*}
\begin{equation*}
\PW\,=\,\left(\, { {1+\delta_{\mu\nu}} \over 2}\, {
{\partial}\over {\partial w_{\mu\nu}} } \,\right),
\end{equation*}
\begin{equation*}
\PWB\,=\,\left(\, { {1+\delta_{\mu\nu}}\over 2} \, {
{\partial}\over {\partial {\overline w}_{\mu\nu} }  } \,\right),
\end{equation*}
$$\PE=\begin{pmatrix} {\partial}\over{\partial \eta_{11}} & \hdots &
 {\partial}\over{\partial \eta_{m1}} \\
\vdots&\ddots&\vdots\\
 {\partial}\over{\partial \eta_{1n}} &\hdots & {\partial}\over
{\partial \eta_{mn}} \end{pmatrix},$$
$$\PEB=\left({\partial}\over{\partial {\overline \eta}_{kl} }\right).$$
%\begin{pmatrix}
%{\partial}\over{\partial {\overline \eta}_{11} }   &
%\hdots&{ {\partial}\over{\partial {\overline \eta}_{m1} }  }\\
%\vdots&\ddots&\vdots\\
%{ {\partial}\over{\partial{\overline \eta}_{1n} }  }&\hdots &
% {\partial}\over{\partial{\overline \eta}_{mn} }  \end{pmatrix}.$$

%\newcommand\ot{\overline\eta}

\noindent {\bf Theorem\ 6 (J.-H. Yang [18], 2005).} The following
metric $d{\tilde s}^2_{n,m;A,B}$ defined by

\end{slide}

%%%%%%%%%%%%%%%%%%%%%%%%%%%%%%%%%%%%%%%%%%%%%%%%%%%%%%%%%%%%%%%%%%%%%%%%%%%%%%%%%%%%%%%%%%%%%%%%%%%%%%%%%%%%%%%%%%%%%%%%%%%%%%%%%%%%%%%
%%%%%%%%%%%%%%%%%%%%%%%%%%%%%%%%%%%%%%%%%%%%%%%%%%%%%%%%%%%%%%%%%%%%%%%%%%%%%%%%%%%%%%%%%%%%%%%%%%%%%%%%%%%%%%%%%%%%%%%%%%%%%%%%%%%%%%%
%%%%%%%%%%%%%%%%%%%%%%%%%%%%%%%%%%%%%%%%%%%%%%%%%%%%%%%%%%%%%%%%%%%%%%%%%%%%%%%%%%%%%%%%%%%%%%%%%%%%%%%%%%%%%%%%%%%%%%%%%%%%%%%%%%%%%%%

\begin{slide}

\newpage
%\vskip 0.2cm
\begin{eqnarray*}
& & {\frac 14}\,d{\tilde s}_{n,m;A,B}^2=\\
& & A\,{\textrm{tr}}\big( (I_n-W\OW)^{-1}dW(I_n-\OW W)^{-1}\bw\big)\,\\
&+&B\, \Big\{ {\textrm{tr}}\big( (I_n-W\OW)^{-1}\,{}^t(d\eta)\,\be\,\big)\\
&+ &{\textrm{tr}}\big(  (\eta\OW-{\overline\eta})(I_n-W\OW)^{-1}dW\\ & &\ \ \ (I_n-\OW W)^{-1}\,{}^t(d\ot)\big)\\
&+&{\textrm{tr}}\big( (\ot W-\eta)(I_n-\OW
W)^{-1}d\OW\\ & &\ \ \ (I_n-W\OW)^{-1}\,{}^t(d\eta)\,\big)    \\
&- &{\textrm{tr}}\big( (I_n-W\OW)^{-1}\,{}^t\eta\,\eta\, (I_n-\OW W)^{-1}\\ &
&\ \ \ \OW
dW (I_n-\OW W)^{-1}d\OW \, \big)\\
&- &{\textrm{tr}}\big( W(I_n-\OW W)^{-1}\,{}^t\ot \,\ot\, (I_n-W\OW )^{-1}\\ &
&
\ \ \ dW (I_n-\OW W)^{-1}d\OW \,\big)\\
&+ &{\textrm{tr}}\big( (I_n-W\OW)^{-1}\,{}^t\eta\,\ot \,(I_n-W\OW)^{-1}\\ & & \ \ \ dW (I_n-\OW W)^{-1} d\OW\,\big)\\
&+ &{\textrm{tr}}\big( (I_n-\OW)^{-1}\,{}^t\ot\,\eta\,\OW\,(I_n-W\OW)^{-1}\\ & & \ \ \ dW (I_n-\OW W)^{-1} d\OW\,\big)\\
&+ &{\textrm{tr}} \big( (I_n-\OW)^{-1}(I_n-W)(I_n-\OW W)^{-1}\\ & & \ \ \
\,{}^t\ot\,\eta\,(I_n-\OW W)^{-1}\, (I_n-\OW)(I_n-W)^{-1}\\ & &
\ \ \ dW(I_n-\OW W)^{-1}d\OW\,\big)\,\Big\}
\end{eqnarray*}

\end{slide}

%%%%%%%%%%%%%%%%%%%%%%%%%%%%%%%%%%%%%%%%%%%%%%%%%%%%%%%%%%%%%%%%%%%%%%%%%%%%%%%%%%%%%%%%%%%%%%%%%%%%%%%%%%%%%%%%%%%%%%%%%%%%%%%%%%%%%%%
%%%%%%%%%%%%%%%%%%%%%%%%%%%%%%%%%%%%%%%%%%%%%%%%%%%%%%%%%%%%%%%%%%%%%%%%%%%%%%%%%%%%%%%%%%%%%%%%%%%%%%%%%%%%%%%%%%%%%%%%%%%%%%%%%%%%%%%
%%%%%%%%%%%%%%%%%%%%%%%%%%%%%%%%%%%%%%%%%%%%%%%%%%%%%%%%%%%%%%%%%%%%%%%%%%%%%%%%%%%%%%%%%%%%%%%%%%%%%%%%%%%%%%%%%%%%%%%%%%%%%%%%%%%%%%%

\begin{slide}

\newpage
\begin{eqnarray*}
%&+ &{\textrm{tr}} \big( (I_n-\OW)^{-1}(I_n-W)(I_n-\OW W)^{-1}\\ & & \ \ \
%\,{}^t\ot\,\eta\,(I_n-\OW W)^{-1}\, (I_n-\OW)(I_n-W)^{-1}\\ & &
%\ \ \ dW(I_n-\OW W)^{-1}d\OW\,\big)\,\Big\}\\
\ &-&B\, {\textrm{tr}} \big((I_n-W\OW)^{-1}(I_n-W)(I_n-\OW)^{-1}\\ & &\ \
\,^t\ot\,\eta\,(I_n-W)^{-1}\,dW (I_n-\OW W)^{-1}d\OW\,\big)
\end{eqnarray*}

\noindent is a Riemannian metric on $\bDnm$ which is invariant
under the action (12) of $G^J_*$.

%\vskip 0.1cm
If $n=m=A=B=1,$ then $d{\tilde s}^2=d{\tilde s}^2_{1,1;1,1}$ is
given by
\begin{eqnarray*}
{\frac 14}\,d{\tilde s}^2&=& { {dW\,d\OW}\over
{(1-|W|^2)^2}}\,+\, { 1 \over {(1-|W|^2)} }\,d\eta\,d\ot\\
&+ &{ {(1+|W|^2)|\eta|^2-\OW \eta^2-W\ot^2}\over {(1-|W|^2)^3} }\,dW\,d\OW\\
&+ &  { {\eta\OW -\ot}\over {(1-|W|^2)^2} }\,dWd\ot\\
& + & {
{\ot W -\eta}\over {(1-|W|^2)^2} }\,d\OW d\eta.
\end{eqnarray*}

\vskip 0.1cm \noindent {\bf Theorem\ 7 (J.-H. Yang [18], 2005).}
The Laplacian ${\tilde \Delta}= {\tilde\Delta}_{n,m;A,B}$ of the
above metric $d{\tilde s}^2_{n,m;A,B}$ is given by

\end{slide}

%%%%%%%%%%%%%%%%%%%%%%%%%%%%%%%%%%%%%%%%%%%%%%%%%%%%%%%%%%%%%%%%%%%%%%%%%%%%%%%%%%%%%%%%%%%%%%%%%%%%%%%%%%%%%%%%%%%%%%%%%%%%%%%%%%%%%%%
%%%%%%%%%%%%%%%%%%%%%%%%%%%%%%%%%%%%%%%%%%%%%%%%%%%%%%%%%%%%%%%%%%%%%%%%%%%%%%%%%%%%%%%%%%%%%%%%%%%%%%%%%%%%%%%%%%%%%%%%%%%%%%%%%%%%%%%
%%%%%%%%%%%%%%%%%%%%%%%%%%%%%%%%%%%%%%%%%%%%%%%%%%%%%%%%%%%%%%%%%%%%%%%%%%%%%%%%%%%%%%%%%%%%%%%%%%%%%%%%%%%%%%%%%%%%%%%%%%%%%%%%%%%%%%%

\begin{slide}

%\newpage
\begin{eqnarray*}
&&{\tilde\Delta}\,=A\,\Biggl\{ \textrm{tr}\left[ (I_n-W\OW)
{}^{{}^{{}^\text{\scriptsize $t$}}}\!\!
\left(
(I_n-W\OW)\PWB\right) \PW\right]\,\\
&+& \textrm{tr} \left[\,{}^t(\eta-\ot\,W)\,
{}^{{}^{{}^\text{\scriptsize $t$}}}\!\!
\left( \PEB\right)
 (I_n-\OW W)\PW  \right]\,\\
&+& \textrm{tr} \left[ (\ot-\eta\,\OW)
{}^{{}^{{}^\text{\scriptsize $t$}}}\!\!
\left(
(I_n-W\OW)\PWB\right)\PE\right] \\
&- &  \textrm{tr}  \left[ \eta \OW (I_n-W\OW)^{-1}
{}^t\eta  \,
{}^{{}^{{}^\text{\scriptsize $t$}}}\!\!
\left(\PEB\right) (I_n-\OW
W)\PE\right]\\
&- &  \textrm{tr}  \left[ \ot W (I_n-\OW W)^{-1}
{}^t\ot\,
{}^{{}^{{}^\text{\scriptsize $t$}}}\!\!
\left(\PEB\right)(I_n-\OW
W)\PE\right]\\
&+&  \textrm{tr}  \left[ \ot (I_n-W\OW)^{-1}\,{}^t\eta\,
{}^{{}^{{}^\text{\scriptsize $t$}}}\!\!
\left( \PEB\right)
(I_n-\OW W)\PE\right]\\
&+& \textrm{tr}  \left[ \eta\OW W (I_n-\OW W)^{-1}
\,{}^t\ot
{}^{{}^{{}^\text{\scriptsize $t$}}}\!\!
\left( \PEB\right)
(I_n-\OW W)\PE\right]\Biggl\}\\
&+& B\cdot  \textrm{tr}   \left[ (I_n-\OW W) \PE\,
{}^{{}^{{}^\text{\scriptsize $t$}}}\!\!
\left( \PEB\right)\right].
\end{eqnarray*}

\end{slide}

%%%%%%%%%%%%%%%%%%%%%%%%%%%%%%%%%%%%%%%%%%%%%%%%%%%%%%%%%%%%%%%%%%%%%%%%%%%%%%%%%%%%%%%%%%%%%%%%%%%%%%%%%%%%%%%%%%%%%%%%%%%%%%%%%%%%%%%
%%%%%%%%%%%%%%%%%%%%%%%%%%%%%%%%%%%%%%%%%%%%%%%%%%%%%%%%%%%%%%%%%%%%%%%%%%%%%%%%%%%%%%%%%%%%%%%%%%%%%%%%%%%%%%%%%%%%%%%%%%%%%%%%%%%%%%%
%%%%%%%%%%%%%%%%%%%%%%%%%%%%%%%%%%%%%%%%%%%%%%%%%%%%%%%%%%%%%%%%%%%%%%%%%%%%%%%%%%%%%%%%%%%%%%%%%%%%%%%%%%%%%%%%%%%%%%%%%%%%%%%%%%%%%%%

\newcommand\ddww{{{\partial^2}\over{\partial W\partial \OW}}}
\newcommand\ddtt{{{\partial^2}\over{\partial \eta\partial \ot}}}
\newcommand\ddwe{{{\partial^2}\over{\partial W\partial \ot}}}
\newcommand\ddew{{{\partial^2}\over{\partial \OW\partial \eta}}}

\begin{slide}

\vskip 0.2cm If $n=m=A=B=1,$ we get
\begin{eqnarray*}
{\tilde\Delta}_{1,1;1,1}&=&\ \ (1-|W|^2)^2 \ddww\,\\ & & +\,(1-|W|^2)\,\ddtt\\
& &+\,(1-|W|^2)(\eta-\ot\,W)\,\ddwe\,\\ & & +\,(1-|W|^2)(\ot-\eta\,\OW)\,\ddew\\
& &-(\OW\,\eta^2+W\ot^2)\,\ddtt\,\\ & & +\,(1+|W|^2)|\eta|^2\,\ddtt.
\end{eqnarray*}

\vskip 0.2cm The main ingredients for the proof of Theorem 6 and
Theorem 7 are the partial Cayley transform (Theorem 5), Theorem 1
and Theorem 2.

\ \ Let ${\Bbb D}(\bDnm)$ be the algebra of all differential
operators $\bDnm$ invariant under the action (12)

\end{slide}

%%%%%%%%%%%%%%%%%%%%%%%%%%%%%%%%%%%%%%%%%%%%%%%%%%%%%%%%%%%%%%%%%%%%%%%%%%%%%%%%%%%%%%%%%%%%%%%%%%%%%%%%%%%%%%%%%%%%%%%%%%%%%%%%%%%%%%%
%%%%%%%%%%%%%%%%%%%%%%%%%%%%%%%%%%%%%%%%%%%%%%%%%%%%%%%%%%%%%%%%%%%%%%%%%%%%%%%%%%%%%%%%%%%%%%%%%%%%%%%%%%%%%%%%%%%%%%%%%%%%%%%%%%%%%%%
%%%%%%%%%%%%%%%%%%%%%%%%%%%%%%%%%%%%%%%%%%%%%%%%%%%%%%%%%%%%%%%%%%%%%%%%%%%%%%%%%%%%%%%%%%%%%%%%%%%%%%%%%%%%%%%%%%%%%%%%%%%%%%%%%%%%%%%

\begin{slide}

of $G^J_*.$ By
Theorem 5, we have the algebra isomorphism
\begin{equation*}
{\Bbb D}(\bDnm)\cong {\Bbb D}(\bHnm).
\end{equation*}

%\newpage
\begin{center}
{\large \bf 6. A fundamental domain for $\G_{n,m}\ba \bHnm$}
\end{center}

\ \ Before we describe a fundamental domain for the Siegel-Jacobi
space, we review the Siegel's fundamental domain for the Siegel
upper half plane. \par\ \ We let
\begin{equation*}
{\mathcal P}_n=\left\{ Y\in\BR^{(n,n)}\,|\ Y=\,{}^tY >0\,\right\}
\end{equation*}
\noindent be an open cone in $\BR^{n(n+1)/2}$. The general linear
group $GL(n,\BR)$ acts on ${\mathcal P}_n$ transitively by
\begin{equation*}
h\circ Y=h\,Y\,{}^th,\quad h\in GL(n,\BR),\ Y\in {\mathcal P}_n.
\end{equation*}

\end{slide}

%%%%%%%%%%%%%%%%%%%%%%%%%%%%%%%%%%%%%%%%%%%%%%%%%%%%%%%%%%%%%%%%%%%%%%%%%%%%%%%%%%%%%%%%%%%%%%%%%%%%%%%%%%%%%%%%%%%%%%%%%%%%%%%%%%%%%%%
%%%%%%%%%%%%%%%%%%%%%%%%%%%%%%%%%%%%%%%%%%%%%%%%%%%%%%%%%%%%%%%%%%%%%%%%%%%%%%%%%%%%%%%%%%%%%%%%%%%%%%%%%%%%%%%%%%%%%%%%%%%%%%%%%%%%%%%
%%%%%%%%%%%%%%%%%%%%%%%%%%%%%%%%%%%%%%%%%%%%%%%%%%%%%%%%%%%%%%%%%%%%%%%%%%%%%%%%%%%%%%%%%%%%%%%%%%%%%%%%%%%%%%%%%%%%%%%%%%%%%%%%%%%%%%%

\begin{slide}

\noindent Thus ${\mathcal P}_n$ is a symmetric space diffeomorphic
to $GL(n,\BR)/O(n).$ We let
\begin{eqnarray*}
GL(n,{\mathbb Z})=\Big\{ h\in GL(n,\BR)\,\Big| \ h\ \textrm{is
integral}\ \Big\}
\end{eqnarray*}
\noindent be the discrete subgroup of $GL(n,\BR)$.

\vskip 0.2cm \ \ The fundamental domain $\Rg$ for $GL(n,\BZ)\ba
{\mathcal P}_n$ which was found by H. Minkowski\,[5] is defined as
a subset of ${\mathcal P}_n$ consisting of $Y=(y_{ij})\in
{\mathcal P}_n$ satisfying the following conditions (M.1)-(M.2)\
(cf.\,[4] p.\,123): \vskip 0.1cm (M.1)\ \ \ $aY\,^ta\geq y_{kk}$\
\ for every $a=(a_i)\in\BZ^n$ in which $a_k,\cdots,a_n$ are
relatively prime for $k=1,2,\cdots,n$. \\ (M.2)\ \ \ \
$y_{k,k+1}\geq 0$ \ for $k=1,\cdots,n-1.$ \par \noindent We say
that a point of $\Rg$ is {\it Minkowski reduced} or simply {\it
M}-{\it reduced}.

\end{slide}

%%%%%%%%%%%%%%%%%%%%%%%%%%%%%%%%%%%%%%%%%%%%%%%%%%%%%%%%%%%%%%%%%%%%%%%%%%%%%%%%%%%%%%%%%%%%%%%%%%%%%%%%%%%%%%%%%%%%%%%%%%%%%%%%%%%%%%%
%%%%%%%%%%%%%%%%%%%%%%%%%%%%%%%%%%%%%%%%%%%%%%%%%%%%%%%%%%%%%%%%%%%%%%%%%%%%%%%%%%%%%%%%%%%%%%%%%%%%%%%%%%%%%%%%%%%%%%%%%%%%%%%%%%%%%%%
%%%%%%%%%%%%%%%%%%%%%%%%%%%%%%%%%%%%%%%%%%%%%%%%%%%%%%%%%%%%%%%%%%%%%%%%%%%%%%%%%%%%%%%%%%%%%%%%%%%%%%%%%%%%%%%%%%%%%%%%%%%%%%%%%%%%%%%

\begin{slide}

\vskip 0.1cm \ \ Siegel\,[8] determined a fundamental domain
${\mathcal F}_n$ for $\G_n\ba \bH_n,$ where $\G_n=Sp(n,\BZ)$ is
the Siegel modular group of degree $n$. We say that $\Om=X+iY\in
\bH_n$ with $X,\,Y$ real is {\it Siegel reduced} or {\it S}-{\it
reduced} if it has the following three properties: \par (S.1)\ \ \
$\det (\text{Im}\,(\g\cdot\Om))\leq \det
(\text{Im}\,(\Om))\quad\text{for\ all}\ \g\in\G_n$; \par (S.2)\ \
$Y=\text{Im}\,\Om$ is M-reduced, that is, $Y\in \Rg\,;$ \par (S.3)
\ \ $|x_{ij}|\leq {\frac 12}\quad \text{for}\ 1\leq i,j\leq n,\
\text{where}\ X=(x_{ij}).$ \vskip 0.21cm \ \ ${\mathcal F}_n$ is
defined as the set of all Siegel reduced points in $\bH_n.$ Using
the highest point method, Siegel [8] proved the following
(F1)-(F3)\,(cf. [4], p.\,169):

\end{slide}

%%%%%%%%%%%%%%%%%%%%%%%%%%%%%%%%%%%%%%%%%%%%%%%%%%%%%%%%%%%%%%%%%%%%%%%%%%%%%%%%%%%%%%%%%%%%%%%%%%%%%%%%%%%%%%%%%%%%%%%%%%%%%%%%%%%%%%%
%%%%%%%%%%%%%%%%%%%%%%%%%%%%%%%%%%%%%%%%%%%%%%%%%%%%%%%%%%%%%%%%%%%%%%%%%%%%%%%%%%%%%%%%%%%%%%%%%%%%%%%%%%%%%%%%%%%%%%%%%%%%%%%%%%%%%%%
%%%%%%%%%%%%%%%%%%%%%%%%%%%%%%%%%%%%%%%%%%%%%%%%%%%%%%%%%%%%%%%%%%%%%%%%%%%%%%%%%%%%%%%%%%%%%%%%%%%%%%%%%%%%%%%%%%%%%%%%%%%%%%%%%%%%%%%

\begin{slide}

\newpage

\vskip 0.1cm (F1)\ \ \ $\G_n\cdot
{\mathcal F}_n=\bH_n,$ i.e., $\bH_n=\cup_{\g\in\G_n}\g\cdot
{\mathcal F}_n.$
\vskip 0.1cm (F2)\ \ ${\mathcal F}_n$ is closed
in $\bH_n.$ \vskip 0.1cm (F3)\ \ ${\mathcal F}_n$ is connected and
the boundary of ${\mathcal F}_n$ consists of a finite number of
hyperplanes.

\newcommand\CCF{{\mathcal F}_n}
\newcommand\CHX{{\mathcal H}_{\xi}}
\newcommand\Ggh{\Gamma_{g,h}}
\newcommand\BHn{\BH_n}
\newcommand\BDn{\BD_n}
\vskip 0.2cm \ \ The metric $ds_{n;1}^2$ induces a metric
$ds_{{\mathcal F}_n}^2$ on ${\mathcal F}_n$. Siegel\,[8] computed
the volume of ${\mathcal F}_n$
\begin{equation*}
\text{vol}\,(\CCF)=2\prod_{k=1}^n\pi^{-k}\G
(k)\zeta(2k),\end{equation*} where $\G (s)$ denotes the Gamma
function and $\zeta (s)$ denotes the Riemann zeta function. For
instance,
\begin{eqnarray*}
& & \text{vol}\,({\mathcal F}_1)={{\pi}\over 3},\quad
\text{vol}\,({\mathcal F}_2)={{\pi^3}\over {270}},\\ &&
\text{vol}\,({\mathcal F}_3)={{\pi^6}\over {127575}},\quad
\text{vol}\,({\mathcal F}_4)={{\pi^{10}}\over {200930625}}.
\end{eqnarray*}

\end{slide}

%%%%%%%%%%%%%%%%%%%%%%%%%%%%%%%%%%%%%%%%%%%%%%%%%%%%%%%%%%%%%%%%%%%%%%%%%%%%%%%%%%%%%%%%%%%%%%%%%%%%%%%%%%%%%%%%%%%%%%%%%%%%%%%%%%%%%%%
%%%%%%%%%%%%%%%%%%%%%%%%%%%%%%%%%%%%%%%%%%%%%%%%%%%%%%%%%%%%%%%%%%%%%%%%%%%%%%%%%%%%%%%%%%%%%%%%%%%%%%%%%%%%%%%%%%%%%%%%%%%%%%%%%%%%%%%
%%%%%%%%%%%%%%%%%%%%%%%%%%%%%%%%%%%%%%%%%%%%%%%%%%%%%%%%%%%%%%%%%%%%%%%%%%%%%%%%%%%%%%%%%%%%%%%%%%%%%%%%%%%%%%%%%%%%%%%%%%%%%%%%%%%%%%%

\begin{slide}

\newpage

\vskip 0.2cm \ \ Let $f_{kl}\,(1\leq k\leq m,\ 1\leq l\leq n)$ be
the $m\times n$ matrix with entry $1$ where the $k$-th row and the
$l$-th column meet, and all other entries $0$. For an element
$\Om\in \BH_n$, we set for brevity
\begin{equation*}
h_{kl}(\Om)=f_{kl}\Om,\quad 1\leq k\leq m,\ 1\leq l\leq
n.\end{equation*}

\ \ For each $\Om\in {\mathcal F}_n,$ we define a subset $P_{\Om}$
of $\BC^{(m,n)}$ by
\begin{eqnarray*}
P_{\Om}=\Big\{ \,\sum_{k=1}^m\sum_{j=1}^n & \lambda_{kl}f_{kl}+
\sum_{k=1}^m\sum_{j=1}^n \mu_{kl}h_{kl}(\Om)\,\Big|\\
& 0\leq \lambda_{kl},\mu_{kl}\leq 1\,\Big\}. \end{eqnarray*}

\noindent \ \ For each $\Om\in {\mathcal F}_n,$ we define the
subset $D_{\Om}$ of $\bHnm$ by
\begin{equation*} D_{\Om}=\big\{\,(\Om,Z)\in\bHnm\,\vert\ Z\in
P_{\Om}\,\big\}.\end{equation*}

\noindent We define
\begin{equation*} \Fgh=\cup_{\Om\in {\mathcal F}_n}D_{\Omega}.\end{equation*}

\end{slide}

%%%%%%%%%%%%%%%%%%%%%%%%%%%%%%%%%%%%%%%%%%%%%%%%%%%%%%%%%%%%%%%%%%%%%%%%%%%%%%%%%%%%%%%%%%%%%%%%%%%%%%%%%%%%%%%%%%%%%%%%%%%%%%%%%%%%%%%
%%%%%%%%%%%%%%%%%%%%%%%%%%%%%%%%%%%%%%%%%%%%%%%%%%%%%%%%%%%%%%%%%%%%%%%%%%%%%%%%%%%%%%%%%%%%%%%%%%%%%%%%%%%%%%%%%%%%%%%%%%%%%%%%%%%%%%%
%%%%%%%%%%%%%%%%%%%%%%%%%%%%%%%%%%%%%%%%%%%%%%%%%%%%%%%%%%%%%%%%%%%%%%%%%%%%%%%%%%%%%%%%%%%%%%%%%%%%%%%%%%%%%%%%%%%%%%%%%%%%%%%%%%%%%%%

\begin{slide}

\vskip 0.3cm \noindent {\bf Theorem\ 8 (J.-H. Yang
[19], 2005).} Let
$$\Gamma_{n,m}=\G_n\ltimes H_{\mathbb Z}^{(n,m)}$$
be the discrete subgroup of $G^J$. Then $\Fgh$ is a fundamental
domain for $\G_{n,m}\ba \bHnm.$ \vskip 0.2cm \noindent
$\textit{Proof.}$ The proof can be found in [19]. \hfill $\square$

%%%%%%%%%%%%%%%%%%%%%%%%%%%%%%%%%%%%%%%%%%%%%%%%%%%%%%%%%%%%%%%%%%%%%%%%%%%%%%%%%%%%%%%%%%%%%%%%%%%%%%%%%
%%%%%%%%%%%%%%%%%%%%%%%%%%%%%%%%%%%%%%%%%%%%%%%%%%%%%%%%%%%%%%%%%%%%%%%%%%%%%%%%%%%%%%%%%%%%%%%%%%%%%%%%%
%%%%%%%%%%%%%%%%%%%%%%%%%%%%%%%%%%%%%%%%%%%%%%%%%%%%%%%%%%%%%%%%%%%%%%%%%%%%%%%%%%%%%%%%%%%%%%%%%%%%%%%%%

**************************************
**************************************
**************************************
%\newpage

\begin{center}
{\large \bf 7. Maass-Jacobi forms}
\end{center}

\noindent {\bf Definition.} For brevity, we set
$\Delta_{n,m}=\Delta_{n,m;1,1}$ (cf. Theorem 2). Let
$$\Gamma_{n,m}=\G_n\ltimes H_{\mathbb Z}^{(n,m)}$$
be the discrete subgroup of $G^J$. A smooth function $f:\bHnm\lrt
\BC$ is called a $\textbf{Maass}$-$\textbf{Jacobi form}$ on
$\bHnm$ if $f$ satisfies the following conditions
(MJ1)-(MJ3)\,:

\end{slide}

%%%%%%%%%%%%%%%%%%%%%%%%%%%%%%%%%%%%%%%%%%%%%%%%%%%%%%%%%%%%%%%%%%%%%%%%%%%%%%%%%%%%%%%%%%%%%%%%%%%%%%%%%%%%%%%%%%%%%%%%%%%%%%%%%%%%%%%
%%%%%%%%%%%%%%%%%%%%%%%%%%%%%%%%%%%%%%%%%%%%%%%%%%%%%%%%%%%%%%%%%%%%%%%%%%%%%%%%%%%%%%%%%%%%%%%%%%%%%%%%%%%%%%%%%%%%%%%%%%%%%%%%%%%%%%%
%%%%%%%%%%%%%%%%%%%%%%%%%%%%%%%%%%%%%%%%%%%%%%%%%%%%%%%%%%%%%%%%%%%%%%%%%%%%%%%%%%%%%%%%%%%%%%%%%%%%%%%%%%%%%%%%%%%%%%%%%%%%%%%%%%%%%%%

\begin{slide}

\vskip 0.1cm (MJ1)\ \ \ $f$ is invariant under
$\G_{n,m}.$\par (MJ2)\ \ \ $f$ is an eigenfunction of
 $\Delta_{n,m}$.\par (MJ3)\ \ \ $f$
has a polynomial growth, that is,\\ there exist a constant $C>0$
and a positive integer $N$ such that
\begin{eqnarray*}
&& |f(X+iY,Z)|\leq C\,|p(Y)|^N \\ && \quad \textrm{as}\ \det
Y\longrightarrow \infty,
\end{eqnarray*}
where $p(Y)$ is a polynomial in $Y=(y_{ij}).$ (cf. See Section 6)

\vskip 0.3cm It is natural to propose the following problems.

\vskip 0.3cm\noindent {\bf {Problem\ C}.} Construct Maass-Jacobi
forms.

\vskip 0.3cm\noindent {\bf {Problem\ D}.} Find all the
eigenfunctions of \\
${}\  \quad {}\qquad {}\qquad\  \Delta_{n,m}.$

\end{slide}

%%%%%%%%%%%%%%%%%%%%%%%%%%%%%%%%%%%%%%%%%%%%%%%%%%%%%%%%%%%%%%%%%%%%%%%%%%%%%%%%%%%%%%%%%%%%%%%%%%%%%%%%%%%%%%%%%%%%%%%%%%%%%%%%%%%%%%%
%%%%%%%%%%%%%%%%%%%%%%%%%%%%%%%%%%%%%%%%%%%%%%%%%%%%%%%%%%%%%%%%%%%%%%%%%%%%%%%%%%%%%%%%%%%%%%%%%%%%%%%%%%%%%%%%%%%%%%%%%%%%%%%%%%%%%%%
%%%%%%%%%%%%%%%%%%%%%%%%%%%%%%%%%%%%%%%%%%%%%%%%%%%%%%%%%%%%%%%%%%%%%%%%%%%%%%%%%%%%%%%%%%%%%%%%%%%%%%%%%%%%%%%%%%%%%%%%%%%%%%%%%%%%%%%

\begin{slide}

\vskip 0.3cm  We consider the simple case $n=m=A=B=1.$ A metric
$ds_{1,1}^2=ds_{1,1;1,1}^2$ on ${\mathbf H}_1\times \BC$ given by
\begin{align*} ds^2_{1,1}\,=\,&{{y\,+\,v^2}\over
{y^3}}\,(\,dx^2\,+\,dy^2\,)\,+\, {\frac 1y}\,(\,du^2\,+\,dv^2\,)\\
&\ \ -\,{{2v}\over {y^2}}\, (\,dx\,du\,+\,dy\,dv\,)\end{align*} is
a $G^J$-invariant K{\"a}hler metric on ${\mathbf H}_1\times \BC$.
Its Laplacian $\Delta_{1,1}$ is given by
\begin{align*} \Delta_{1,1}\,=\,& \,y^2\,\left(\,\ddx\,+\,\ddy\,\right)\,\\
&+\, (\,y\,+\,v^2\,)\,\left(\,\ddu\,+\,\ddv\,\right)\\ &\ \
+\,2\,y\,v\,\left(\,\pxu\,+\,\pyv\,\right).\end{align*}

\vskip 0.2cm We provide some examples of eigenfunctions of
$\Delta_{1,1}$.

\end{slide}

%%%%%%%%%%%%%%%%%%%%%%%%%%%%%%%%%%%%%%%%%%%%%%%%%%%%%%%%%%%%%%%%%%%%%%%%%%%%%%%%%%%%%%%%%%%%%%%%%%%%%%%%%%%%%%%%%%%%%%%%%%%%%%%%%%%%%%%
%%%%%%%%%%%%%%%%%%%%%%%%%%%%%%%%%%%%%%%%%%%%%%%%%%%%%%%%%%%%%%%%%%%%%%%%%%%%%%%%%%%%%%%%%%%%%%%%%%%%%%%%%%%%%%%%%%%%%%%%%%%%%%%%%%%%%%%
%%%%%%%%%%%%%%%%%%%%%%%%%%%%%%%%%%%%%%%%%%%%%%%%%%%%%%%%%%%%%%%%%%%%%%%%%%%%%%%%%%%%%%%%%%%%%%%%%%%%%%%%%%%%%%%%%%%%%%%%%%%%%%%%%%%%%%%

\begin{slide}

\vskip 0.2cm (1) $h(x,y)=y^{1\over
2}K_{s-{\frac12}}(2\pi |a|y)\,e^{2\pi iax} \ (s\in \BC,$
$a\not=0\,)$ with eigenvalue $s(s-1).$ Here
$$K_s(z):={\frac12}\int^{\infty}_0 \exp\left\{-{z\over
2}(t+t^{-1})\right\}\,t^{s-1}\,dt,$$ \indent \ \ \ where
$\mathrm{Re}\,z
> 0.$ \par (2) $y^s,\ y^s x,\ y^s u\ (s\in\BC)$ with eigenvalue
$s(s-1).$
\par
 (3) $y^s v,\ y^s uv,\ y^s xv$ with eigenvalue $s(s+1).$
\par
(4) $x,\,y,\,u,\,v,\,xv,\,uv$ with eigenvalue $0$.
\par
(5) All Maass wave forms.

\vskip 0.5cm\noindent {\bf 7.1. Eisenstein Series} \par\ \  Let
$$\Gamma_{1,1}^{\infty}=\left\{ \left( \begin{pmatrix} \pm 1 & m\\
0 & \pm 1\end{pmatrix},(0,n,\kappa)\right)\,\Big|\ m,n,\kappa\in
\BZ\,\right\}$$

\end{slide}

%%%%%%%%%%%%%%%%%%%%%%%%%%%%%%%%%%%%%%%%%%%%%%%%%%%%%%%%%%%%%%%%%%%%%%%%%%%%%%%%%%%%%%%%%%%%%%%%%%%%%%%%%%%%%%%%%%%%%%%%%%%%%%%%%%%%%%%
%%%%%%%%%%%%%%%%%%%%%%%%%%%%%%%%%%%%%%%%%%%%%%%%%%%%%%%%%%%%%%%%%%%%%%%%%%%%%%%%%%%%%%%%%%%%%%%%%%%%%%%%%%%%%%%%%%%%%%%%%%%%%%%%%%%%%%%
%%%%%%%%%%%%%%%%%%%%%%%%%%%%%%%%%%%%%%%%%%%%%%%%%%%%%%%%%%%%%%%%%%%%%%%%%%%%%%%%%%%%%%%%%%%%%%%%%%%%%%%%%%%%%%%%%%%%%%%%%%%%%%%%%%%%%%%

\begin{slide}

be the subgroup of $\Gamma_{1,1}=SL_2(\BZ)\ltimes
H_{\BZ}^{(1,1)}.$
\par
\ \ \ For $\gamma=\left( \begin{pmatrix} a & b\\
c & d\end{pmatrix},(\lambda,\mu,\kappa)\right)\in \Gamma_{1,1},$
we put \\ $(\tau_{\gamma},z_{\gamma})=\gamma\cdot (\tau,z).$ That
is,
\begin{align*}
\tau_{\gamma}&=(a\tau+b)(c\tau+d)^{-1},\\ z_{\gamma}&=(z+\lambda
\tau+\nu)(c\tau+d)^{-1}.\end{align*} We note that if
$\gamma\in\Gamma_{1,1},$
$$\mathrm{Im}\,\tau_{\gamma}=\mathrm{Im}\,\tau,\ \
\mathrm{Im}\,z_{\gamma}=\mathrm{Im}\,z$$ if and only if
$\gamma\in\Gamma_{1,1}^{\infty}.$ For $s\in\BC,$ we define an
Eisenstein series formally by
$$E_s(\tau,z)=\sum_{\gamma\in \Gamma_{1,1}^{\infty}\backslash
\Gamma_{1,1}}(\mathrm{Im}\,\tau_{\gamma})^s\cdot
\mathrm{Im}\,z_{\gamma}.$$ Then $E_s(\tau,z)$ satisfies formally
$$E_s(\gamma\cdot(\tau,z))=E_s(\tau,z),\ \ \gamma\in
\Gamma_{1,1}$$ and $$\Delta E_s(\tau,z)=s(s+1)E_s(\tau,z).$$
%But
%the series does not converge.

\end{slide}

%%%%%%%%%%%%%%%%%%%%%%%%%%%%%%%%%%%%%%%%%%%%%%%%%%%%%%%%%%%%%%%%%%%%%%%%%%%%%%%%%%%%%%%%%%%%%%%%%%%%%%%%%%%%%%%%%%%%%%%%%%%%%%%%%%%%%%%
%%%%%%%%%%%%%%%%%%%%%%%%%%%%%%%%%%%%%%%%%%%%%%%%%%%%%%%%%%%%%%%%%%%%%%%%%%%%%%%%%%%%%%%%%%%%%%%%%%%%%%%%%%%%%%%%%%%%%%%%%%%%%%%%%%%%%%%
%%%%%%%%%%%%%%%%%%%%%%%%%%%%%%%%%%%%%%%%%%%%%%%%%%%%%%%%%%%%%%%%%%%%%%%%%%%%%%%%%%%%%%%%%%%%%%%%%%%%%%%%%%%%%%%%%%%%%%%%%%%%%%%%%%%%%%%

\begin{slide}

\vskip 0.5cm\noindent {\bf 7.2. Fourier Expansion of Maass-Jacobi
\\ Form}
\newcommand{\bHC}{{\mathbf H}_1\times \BC}

\ \ We let $f:{\mathbf H}_1\times \BC\lrt \BC$ be a Maass-Jacobi
form with $\Delta f=\lambda f.$ Then $f$ satisfies the following
invariance relations
$$f(\tau+n,\,z)\,=\,f(\tau,z)\ \ \ \mathrm{for\ all}\ n\in
\BZ$$ and
$$f(\tau,\,z\,+\,n_1\tau\,+\,n_2)\,=\,f(\tau,z)$$
$\mathrm{for\ all}\ n_1,\,n_2\in \BZ.$ Therefore $f$ is a smooth
function on $\bHC$ which is periodic in $x$ and $u$ with period
$1.$ So $f$ has the following Fourier series
$$f(\tau,z)\,=\,\sum_{n\in \BZ}\sum_{r\in \BZ}\,c_{n,r}(y,v)\,
e^{2\pi i(nx+ru)}.$$ For two fixed integers $n$ and $r$, we have
to calculate the function $c_{n,r}(y,v).$ For brevity, we

\end{slide}

%%%%%%%%%%%%%%%%%%%%%%%%%%%%%%%%%%%%%%%%%%%%%%%%%%%%%%%%%%%%%%%%%%%%%%%%%%%%%%%%%%%%%%%%%%%%%%%%%%%%%%%%%%%%%%%%%%%%%%%%%%%%%%%%%%%%%%%
%%%%%%%%%%%%%%%%%%%%%%%%%%%%%%%%%%%%%%%%%%%%%%%%%%%%%%%%%%%%%%%%%%%%%%%%%%%%%%%%%%%%%%%%%%%%%%%%%%%%%%%%%%%%%%%%%%%%%%%%%%%%%%%%%%%%%%%
%%%%%%%%%%%%%%%%%%%%%%%%%%%%%%%%%%%%%%%%%%%%%%%%%%%%%%%%%%%%%%%%%%%%%%%%%%%%%%%%%%%%%%%%%%%%%%%%%%%%%%%%%%%%%%%%%%%%%%%%%%%%%%%%%%%%%%%

\begin{slide}

put
$F(y,v)=\,c_{n,r}(y,v).$ Then $F$ satisfies the following
differential equation
\begin{align*}
&\left[y^2\,\ddy\,+\,(y+v^2)\,\ddv\,+\,2yv\, \pyv\,\right]\,F\\ =
& \ \left\{\,(ay+bv)^2\,+\,b^2y\,+\,\lambda\,\right\}
\,F.\end{align*} Here $a=2\pi n$ and $b=2\pi r$ are constant. We
note that the function $u(y)=y^{\frac 12} K_{s-{\frac 12}}(2\pi
\vert n\vert y)$ satisfies the above differential equation with
$\lambda=s(s-1).$ Here $K_s(z)$ is the $K$-Bessel function before.

{\bf \underline{Problem}\,:} Find the solutions of the above
differential equations explicitly.

{\bf \underline{Problem}\,:} Develop a Fourier expansion of a Maass-Jacobi form in terms
of the Whittaker functions.

**************************************
**************************************

\end{slide}

%%%%%%%%%%%%%%%%%%%%%%%%%%%%%%%%%%%%%%%%%%%%%%%%%%%%%%%%%%%%%%%%%%%%%%%%%%%%%%%%%%%%%%%%%%%%%%%%%%%%%%%%%%%%%%%%%%%%%%%%%%%%%%%%%%%%%%%
%%%%%%%%%%%%%%%%%%%%%%%%%%%%%%%%%%%%%%%%%%%%%%%%%%%%%%%%%%%%%%%%%%%%%%%%%%%%%%%%%%%%%%%%%%%%%%%%%%%%%%%%%%%%%%%%%%%%%%%%%%%%%%%%%%%%%%%
%%%%%%%%%%%%%%%%%%%%%%%%%%%%%%%%%%%%%%%%%%%%%%%%%%%%%%%%%%%%%%%%%%%%%%%%%%%%%%%%%%%%%%%%%%%%%%%%%%%%%%%%%%%%%%%%%%%%%%%%%%%%%%%%%%%%%%%

\begin{slide}

\newpage

%Discuss the existence and uniqueness of the
%Whittaker model (via an integral transform).

\begin{center}
{\large \bf 8. Spectral theory of $\Delta_{n,m;A,B}$}
\end{center}
\par
$\textbf{Problem :}$ Develop the spectral theory of $\Delta_{n,m}$
on ${\mathcal F}_{n,m}$.

$ \textbf{Step I.}$  Spectral Theory of $\Delta_{\Om}$ on
$A_{\Om}$

$ \textbf{Step II.}$ Spectral Theory of $\Delta_n$ on ${\mathcal
F}_n$ (Hard at this moment)

$ \textbf{Step III.}$ Mixed Spectral Theory

$ \textbf{Step IV.}$ Combine Step I-III and more advanced works to
develop the Spectral Theory of $\Delta_{n,m}$ on ${\mathcal
F}_{n,m}$.\\
$ \textbf{[Very Complicated and Hard at}$ $ \textbf{this moment]}$

\end{slide}

%%%%%%%%%%%%%%%%%%%%%%%%%%%%%%%%%%%%%%%%%%%%%%%%%%%%%%%%%%%%%%%%%%%%%%%%%%%%%%%%%%%%%%%%%%%%%%%%%%%%%%%%%%%%%%%%%%%%%%%%%%%%%%%%%%%%%%%
%%%%%%%%%%%%%%%%%%%%%%%%%%%%%%%%%%%%%%%%%%%%%%%%%%%%%%%%%%%%%%%%%%%%%%%%%%%%%%%%%%%%%%%%%%%%%%%%%%%%%%%%%%%%%%%%%%%%%%%%%%%%%%%%%%%%%%%
%%%%%%%%%%%%%%%%%%%%%%%%%%%%%%%%%%%%%%%%%%%%%%%%%%%%%%%%%%%%%%%%%%%%%%%%%%%%%%%%%%%%%%%%%%%%%%%%%%%%%%%%%%%%%%%%%%%%%%%%%%%%%%%%%%%%%%%

\begin{slide}

\newpage
\ \ \ I will explain Step I-IV in more detail.

$ \textbf{[Step I]}$ For a fixed element $\Om\in \bH_n,$ we set
\begin{equation*}
L_{\Om}=\BZ^{(m,n)}+\BZ^{(m,n)}\Om\end{equation*} Then $L_{\Om}$
is a lattice in $\BC^{(m,n)}$ and the period matrix
$\Om_*=(I_n,\Om)$ satisfies the Riemann conditions (RC.1) and
(RC.2)\,: \vskip 0.1cm (RC.1) \ \ \ $\Om_*J_n\,\Om_*^T=0\,$;
\vskip 0.1cm (RC.2) \ \ \ $-{1 \over
{i}}\Om_*J_n\,{\overline{\Om}}_*^T
>0$.

Thus the complex torus $A_{\Om}=\BC^{(m,n)}/L_{\Omega}$ is an
abelian variety. For more details on $A_{\Om}$, we refer to [6].

We write $\Om=X+iY$ of $\bH_n$ with $X=\text{Re}\,\Om$ and
$Y=\text{Im}\, \Om.$ For a pair $(A,B)$ with $A,B\in\BZ^{(m,n)},$

\end{slide}

%%%%%%%%%%%%%%%%%%%%%%%%%%%%%%%%%%%%%%%%%%%%%%%%%%%%%%%%%%%%%%%%%%%%%%%%%%%%%%%%%%%%%%%%%%%%%%%%%%%%%%%%%%%%%%%%%%%%%%%%%%%%%%%%%%%%%%%
%%%%%%%%%%%%%%%%%%%%%%%%%%%%%%%%%%%%%%%%%%%%%%%%%%%%%%%%%%%%%%%%%%%%%%%%%%%%%%%%%%%%%%%%%%%%%%%%%%%%%%%%%%%%%%%%%%%%%%%%%%%%%%%%%%%%%%%
%%%%%%%%%%%%%%%%%%%%%%%%%%%%%%%%%%%%%%%%%%%%%%%%%%%%%%%%%%%%%%%%%%%%%%%%%%%%%%%%%%%%%%%%%%%%%%%%%%%%%%%%%%%%%%%%%%%%%%%%%%%%%%%%%%%%%%%

\begin{slide}

we define the function $E_{\Om;A,B}:\Cmn\lrt \BC$ by
\begin{equation*}
E_{\Om;A,B}(Z)=e^{2\pi i\left( \textrm{tr}\,(A^TU\,)+\,\textrm{tr}\,
\big((B-AX)Y^{-1}V^T\big)\right)},\end{equation*} where $Z=U+iV$
is a variable in $\Cmn$ with real $U,V$.

\newcommand\AO{A_{\Omega}}

$ \textbf{Theorem\,:}$ The set $\left\{\,E_{\Om;A,B}\,|\
A,B\in\BZ^{(m,n)}\,\right\}$ is a complete orthonormal basis for
$L^2(\AO)$. Moreover we have the following spectral decomposition
of $\Delta_{\Om}$:
$$L^2(\AO)=\oplus_{A,B\in \BZ^{(m,n)}}\BC\cdot E_{\Om;A,B}.$$

$ \textbf{[Step II]}$  The inner product $(\ ,\ )$ on
$L^2({\mathcal F}_n)$ is defined by
$$(f,g)=\int_{{\mathcal
F}_n} f(\Om)\,\overline{g(\Omega)}\,\,\,{{[dX]\wedge [dY]}\over
{(\det Y)^{n+1}}}.$$

\ \ \ $L^2({\mathcal F}_n)$ is decomposed as follows\,:
\begin{equation*}
L^2({\mathcal F}_n)=L^2_{ \texttt{cusp}}({\mathcal F}_n)\oplus
L^2_{ \texttt{res}}({\mathcal F}_n)\oplus L^2_{
\texttt{cont}}({\mathcal F}_n)
\end{equation*}

\end{slide}

%%%%%%%%%%%%%%%%%%%%%%%%%%%%%%%%%%%%%%%%%%%%%%%%%%%%%%%%%%%%%%%%%%%%%%%%%%%%%%%%%%%%%%%%%%%%%%%%%%%%%%%%%%%%%%%%%%%%%%%%%%%%%%%%%%%%%%%
%%%%%%%%%%%%%%%%%%%%%%%%%%%%%%%%%%%%%%%%%%%%%%%%%%%%%%%%%%%%%%%%%%%%%%%%%%%%%%%%%%%%%%%%%%%%%%%%%%%%%%%%%%%%%%%%%%%%%%%%%%%%%%%%%%%%%%%
%%%%%%%%%%%%%%%%%%%%%%%%%%%%%%%%%%%%%%%%%%%%%%%%%%%%%%%%%%%%%%%%%%%%%%%%%%%%%%%%%%%%%%%%%%%%%%%%%%%%%%%%%%%%%%%%%%%%%%%%%%%%%%%%%%%%%%%

\begin{slide}

The continuous part $L^2_{ \texttt{cont}}({\mathcal F}_n)$
can be understood by the theory of $ \textbf{Eisenstein series}$
developed by Alte Selberg and Robert Langlands. Also the residual
part $L^2_{ \texttt{res}}({\mathcal F}_n)$ can be understood. But
the cuspidal part $L^2_{ \texttt{cusp}}({\mathcal F}_n)$ has not
been well developed yet. We have little knowledge of $
\textbf{cusp forms}$.

For instance, if $n=1$, then every element $f$ in $L^2({\mathcal
F}_1)$ is decomposed into
\begin{equation*}
f=\sum_{n=0}^{\infty} (f,g_n)g_n +{1 \over {4\pi i}}\int_{
\textrm{Re}\,s={\frac 12}} (f,E_s)\,E_s\,ds
\end{equation*}
Here $g_0=\sqrt{\frac 3\pi}\,,\ \left\{ g_n\,|\ n\geq 1\,\right\}$
is an orthonormal basis consisting of $ \textbf{cusp Maass
forms}$. The Eisenstein series $E_s\,(s\in \BC)$ is defined by
$$E_s(\Om)=\sum_{\g\in \G_1(\infty)\backslash \G_1}\left(
\textrm{Im}\,(\g\cdot\Om)\right)^s,\ \ \Om\in\BH_1$$ Here
$\G_1=Sp(1,\BZ)=SL(2,\BZ)$ and
$$\G_1(\infty)=\left\{\,\g\in \G_1\,|\
\g\cdot\infty=\infty\,\right\}.$$

\end{slide}

%%%%%%%%%%%%%%%%%%%%%%%%%%%%%%%%%%%%%%%%%%%%%%%%%%%%%%%%%%%%%%%%%%%%%%%%%%%%%%%%%%%%%%%%%%%%%%%%%%%%%%%%%%%%%%%%%%%%%%%%%%%%%%%%%%%%%%%
%%%%%%%%%%%%%%%%%%%%%%%%%%%%%%%%%%%%%%%%%%%%%%%%%%%%%%%%%%%%%%%%%%%%%%%%%%%%%%%%%%%%%%%%%%%%%%%%%%%%%%%%%%%%%%%%%%%%%%%%%%%%%%%%%%%%%%%
%%%%%%%%%%%%%%%%%%%%%%%%%%%%%%%%%%%%%%%%%%%%%%%%%%%%%%%%%%%%%%%%%%%%%%%%%%%%%%%%%%%%%%%%%%%%%%%%%%%%%%%%%%%%%%%%%%%%%%%%%%%%%%%%%%%%%%%

\begin{slide}

$ \textbf{[Step III-IV]}$ The inner product $(\ ,\ )_{n,m}$ on
$L^2({\mathcal F}_{n,m})$ is defined by
\begin{eqnarray*} && (f,g)_{n,m}\\
&=&\int_{{\mathcal F}_{n,m}}
f(\Om,Z)\,\overline{g(\Omega,Z)}\,\,\,{{[dX] [dY] [dU] [dV]}\over
{(\det Y)^{n+m+1}}}.\end{eqnarray*}

\ \ \ $L^2({\mathcal F}_{n,m})$ is decomposed into
\begin{equation*}
L^2({\mathcal F}_{n,m})=L^2_{ \texttt{cusp}}\oplus L^2_{
\texttt{res}}\oplus L^2_{ \texttt{cont}}
\end{equation*}
The continuous part $L^2_{ \texttt{cont}}$ can be understood by
the theory of Eisenstein series with some more work. But the
cuspidal part $L^2_{ \texttt{cusp}}$ has not been developed yet.
This part is closely related to the theory of $
\textbf{Maass-Jacobi cusp forms}$.

We have the following natural question\,:

$ \textbf{Problem.}$ Develop the theory of Maass-Jacobi forms
(e.g., Hecke theory of Maass-Jacobi forms, Whittaker functions etc).

\end{slide}

%%%%%%%%%%%%%%%%%%%%%%%%%%%%%%%%%%%%%%%%%%%%%%%%%%%%%%%%%%%%%%%%%%%%%%%%%%%%%%%%%%%%%%%%%%%%%%%%%%%%%%%%%%%%%%%%%%%%%%%%%%%%%%%%%%%%%%%
%%%%%%%%%%%%%%%%%%%%%%%%%%%%%%%%%%%%%%%%%%%%%%%%%%%%%%%%%%%%%%%%%%%%%%%%%%%%%%%%%%%%%%%%%%%%%%%%%%%%%%%%%%%%%%%%%%%%%%%%%%%%%%%%%%%%%%%
%%%%%%%%%%%%%%%%%%%%%%%%%%%%%%%%%%%%%%%%%%%%%%%%%%%%%%%%%%%%%%%%%%%%%%%%%%%%%%%%%%%%%%%%%%%%%%%%%%%%%%%%%%%%%%%%%%%%%%%%%%%%%%%%%%%%%%%

\begin{slide}

\begin{center}
{\large \bf 9. Decomposition of the regular representation of
$G^J$}
\end{center}

\ \ \ It is very important to decompose the $  \textbf{regular}$
representation of $G^J$ on $L^2\big(\G_{n,m}\backslash G^J\big)$
into irreducible (unitary) representations. Here
$$\G_{n,m}=Sp(n,\BZ)\ltimes H_{\BZ}^{(n,m)}.$$
\ \ For brevity, we put
$$L^2=L^2\big(\G_{n,m}\backslash G^J\big).$$
\ \ Then the regular representation of $G^J$ is decomposed into
\begin{equation*}
L^2=L^2_d\oplus L^2_c,
\end{equation*}
\noindent where $L^2_d$ is the discrete part of $L^2$ and $L^2_c$
is the continuous part of $L^2.$ The continuous part of $L^2$ can
be understood by the Langlands' theory of Eisenstein series with
some more work. We decompose $L^2_d$ as
\begin{equation*}
L^2_d=\sum_{\pi}m_{\pi}\pi.
\end{equation*}

\end{slide}

%%%%%%%%%%%%%%%%%%%%%%%%%%%%%%%%%%%%%%%%%%%%%%%%%%%%%%%%%%%%%%%%%%%%%%%%%%%%%%%%%%%%%%%%%%%%%%%%%%%%%%%%%%%%%%%%%%%%%%%%%%%%%%%%%%%%%%%
%%%%%%%%%%%%%%%%%%%%%%%%%%%%%%%%%%%%%%%%%%%%%%%%%%%%%%%%%%%%%%%%%%%%%%%%%%%%%%%%%%%%%%%%%%%%%%%%%%%%%%%%%%%%%%%%%%%%%%%%%%%%%%%%%%%%%%%
%%%%%%%%%%%%%%%%%%%%%%%%%%%%%%%%%%%%%%%%%%%%%%%%%%%%%%%%%%%%%%%%%%%%%%%%%%%%%%%%%%%%%%%%%%%%%%%%%%%%%%%%%%%%%%%%%%%%%%%%%%%%%%%%%%%%%%%

\newpage
\begin{slide}

\begin{center}
{\large \bf 10. Open Problems}
\end{center}

\vskip 0.2cm We list the problems to be investigated in the
future.

$ \textbf{Problem 1.}$ Find explicit algebraically independent generators of ${\mathbb D}
({\mathbf H}_{n,m}).$

$ \textbf{Problem 2.}$ Find explicit algebraically independent generators of
$\textrm{Pol}_{m,n}^K=\textrm{Pol}(T_{n,m})^K.$ Here $K=U(n).$
Decompose the representation $\rho$ of $K$ or $K_\BC=GL(n,\BC)$ on
$\textrm{Pol}(T_{n,m})$ explicitly. More precisely if
$$ \rho=\sum_{\sigma\in {\widehat K}} m_\sigma \,\sigma$$
we want to know the multiplicity $m_\sigma$. I think that the representation is not
multiplicity free.

$ [\textbf{Remark]:}$ For a positive integer $r$, we let \\
$\textrm{Pol}_{[r]}(T_n)$ denote the subspace
of $\textrm{Pol}(T_n)$ consisting of homogeneous polynomial functions

\end{slide}

%%%%%%%%%%%%%%%%%%%%%%%%%%%%%%%%%%%%%%%%%%%%%%%%%%%%%%%%%%%%%%%%%%%%%%%%%%%%%%%%%%%%%%%%%%%%%%%%%%%%%%%%%%%%%%%%%%%%%%%%%%%%%%%%%%%%%%%
%%%%%%%%%%%%%%%%%%%%%%%%%%%%%%%%%%%%%%%%%%%%%%%%%%%%%%%%%%%%%%%%%%%%%%%%%%%%%%%%%%%%%%%%%%%%%%%%%%%%%%%%%%%%%%%%%%%%%%%%%%%%%%%%%%%%%%%
%%%%%%%%%%%%%%%%%%%%%%%%%%%%%%%%%%%%%%%%%%%%%%%%%%%%%%%%%%%%%%%%%%%%%%%%%%%%%%%%%%%%%%%%%%%%%%%%%%%%%%%%%%%%%%%%%%%%%%%%%%%%%%%%%%%%%%%

\begin{slide}

on $T_n$ of degree $r$.
The action of $K$ or $K_\BC$ on $\textrm{Pol}_{[r]}(T_n)$ is multiplicity-free (cf. L. Hua,
W. Schmid, G. Shimura et al).

$ \textbf{Problem 3.}$ Let $(\Omega_1,Z_1)$ and $(\Omega_2,Z_2)$ be two given points in
${\mathbf H}_{n,m}$. Express the distance between $(\Omega_1,Z_1)$ and $(\Omega_2,Z_2)$
for the metric $ds^2_{n,m;A,B}$ explicitly.

$ \textbf{Problem 4.}$ Compute the multiplicity $m_{\pi}$ in
$L^2_d=\sum_{\pi}m_{\pi}\pi$ in Section 9. Investigate
the unitary dual of $G^J.$

$ [\textbf{Remark]:}$ The unitary dual of $Sp(n,\BR)$ is not known
for $n\geq 3.$

$ \textbf{Problem 5.}$ Investigate the Schr{\"o}dinger-Weil
representations of $G^J$ in detail.

$ \textbf{Problem 6.}$ Develop the theory of the orbit method for
$G^J$.

\end{slide}

%%%%%%%%%%%%%%%%%%%%%%%%%%%%%%%%%%%%%%%%%%%%%%%%%%%%%%%%%%%%%%%%%%%%%%%%%%%%%%%%%%%%%%%%%%%%%%%%%%%%%%%%%%%%%%%%%%%%%%%%%%%%%%%%%%%%%%%
%%%%%%%%%%%%%%%%%%%%%%%%%%%%%%%%%%%%%%%%%%%%%%%%%%%%%%%%%%%%%%%%%%%%%%%%%%%%%%%%%%%%%%%%%%%%%%%%%%%%%%%%%%%%%%%%%%%%%%%%%%%%%%%%%%%%%%%
%%%%%%%%%%%%%%%%%%%%%%%%%%%%%%%%%%%%%%%%%%%%%%%%%%%%%%%%%%%%%%%%%%%%%%%%%%%%%%%%%%%%%%%%%%%%%%%%%%%%%%%%%%%%%%%%%%%%%%%%%%%%%%%%%%%%%%%

\begin{slide}

$ \textbf{Problem 7.}$ Find the trace formula for $G^J$
with respect to $\G_{n,m}$.

$ \textbf{Problem 8.}$ Find Weyl's law for $G^J.$ Discuss the existence of
nonzero Maass-Jacobi cusp forms.

$ \textbf{Problem 9.}$ Describe the Fourier transform, the inversion formula, the Plancherel formula
and the spherical transform explicitly.

$ \textbf{Problem 10.}$ Discuss the existence and uniqueness of the Whittaker model (e.g., via
an integral transform). In the case $n=m=1$, R. Berndt and R. Schmidt gave two methods to obtain the Whittaker
models (1) by the infinitesimal method and the the method of differential operators, and (2) via an integral
transform [cf.\,Progress in Math. Vol. 163 (1998], pp.\,63-73]. 

\end{slide}

\begin{slide}

\newpage
\begin{center}
{\Large \bf V. References}
\end{center}

\vskip 1cm

[1] $ \textbf{R. Langlands}$, {\em On the Functional Equations
Satisfied by Eisenstein Series}, Lecture Notes in Math. 544,
Springer-Verlag, Berlin and New York (1976).

[2] $ \textbf{H. Maass}$, {\em {\"U}ber eine neue Art von
nichtanalytischen automorphen Funktionen und die Bestimmung
Dirichletscher Reihen durch Funtionalgleichungen,} Math. Ann. {\bf
121} (1949), 141-183.

[3] $ \textbf{H. Maass}$, {\em Die Differentialgleichungen in der
Theorie der Siegelschen Modulfunktionen}, Math. Ann. {\bf 26}
(1953), 44--68.

[4] $ \textbf{H. Maass}$, {\em Siegel modular functions and
Dirichlet series}, Lecture Notes in Math. 216, Springer-Verlag,
Berlin and New York (1971).

\end{slide}

%%%%%%%%%%%%%%%%%%%%%%%%%%%%%%%%%%%%%%%%%%%%%%%%%%%%%%%%%%%%%%%%%%%%%%%%%%%%%%%%%%%%%%%%%%%%%%%%%%%%%%%%%%%%%%%%%%%%%%%%%%%%%%%%%%%%%%%
%%%%%%%%%%%%%%%%%%%%%%%%%%%%%%%%%%%%%%%%%%%%%%%%%%%%%%%%%%%%%%%%%%%%%%%%%%%%%%%%%%%%%%%%%%%%%%%%%%%%%%%%%%%%%%%%%%%%%%%%%%%%%%%%%%%%%%%
%%%%%%%%%%%%%%%%%%%%%%%%%%%%%%%%%%%%%%%%%%%%%%%%%%%%%%%%%%%%%%%%%%%%%%%%%%%%%%%%%%%%%%%%%%%%%%%%%%%%%%%%%%%%%%%%%%%%%%%%%%%%%%%%%%%%%%%

\begin{slide}

[5] $ \textbf{H.~Minkowski}$, \emph{Gesammelte Abhandlungen,}
Chelsea, New York (1967).

[6] $ \textbf{D.~Mumford}$, \emph{Tata Lectures on Theta I,}
Progress in Math. {\bf 28}, Boston-Basel-Stuttgart (1983).

[7] $ \textbf{A. Selberg}$,  \emph{Harmonic analysis and
discontinuous groups in weakly symmetric Riemannian spaces with
applications to Dirichlet series}, J. Indian Math. Soc. B.  {\bf
20} (1956), 47-87.

[8] $ \textbf{C. L. Siegel}$, \emph{Symplectic Geometry,} Amer. J.
Math. {\bf 65} (1943), 1-86; Academic Press, New York and London
(1964); Gesammelte Abhandlungen, no.\,\,{\bf 41,\ vol.\,II},
Springer-Verlag (1966), 274-359.

[9] $ \textbf{C. L. Siegel}$, \emph{Topics in Complex Function
Theory}, Wiley-Interscience, New York, vol. III (1973).

\end{slide}

%%%%%%%%%%%%%%%%%%%%%%%%%%%%%%%%%%%%%%%%%%%%%%%%%%%%%%%%%%%%%%%%%%%%%%%%%%%%%%%%%%%%%%%%%%%%%%%%%%%%%%%%%%%%%%%%%%%%%%%%%%%%%%%%%%%%%%%
%%%%%%%%%%%%%%%%%%%%%%%%%%%%%%%%%%%%%%%%%%%%%%%%%%%%%%%%%%%%%%%%%%%%%%%%%%%%%%%%%%%%%%%%%%%%%%%%%%%%%%%%%%%%%%%%%%%%%%%%%%%%%%%%%%%%%%%
%%%%%%%%%%%%%%%%%%%%%%%%%%%%%%%%%%%%%%%%%%%%%%%%%%%%%%%%%%%%%%%%%%%%%%%%%%%%%%%%%%%%%%%%%%%%%%%%%%%%%%%%%%%%%%%%%%%%%%%%%%%%%%%%%%%%%%%

\begin{slide}

%[10] $ \textbf{A. Terras}$, \emph{Harmonic Analysis on Symmetric
%Spaces and Applications II.}, Springer-Verlag (1988).

[10] $ \textbf{A. Weil}$, \emph{Sur certains groupes d'operateurs
unitaires (French)}, Acta Math. $ \textbf{111}$ (1964), 143--211.

[11] $ \textbf{J.-H. Yang}$, \emph{The Siegel-Jacobi operator,}
Abh. Math. Sem. Univ. Hamburg {\bf 59} (1993), 135--146.

[12] $ \textbf{J.-H. Yang}$, \emph{Singular Jacobi Forms,} Trans.
Amer. Math. Soc. {\bf 347\, (6)} (1995), 2041--2049.

[13] $ \textbf{J.-H. Yang}$, {\em Construction of vector valued
modular forms from Jacobi forms,} Canadian J. of Math. {\bf 47\,
(6)} (1995), 1329--1339.

[14] $ \textbf{J.-H. Yang}$, \emph{A Geometrical Theory of Jacobi
Forms of Higher Degree,} Proceedings of Symposium on Hodge Theory
and Algebraic Geometry (Editor: Tadao Oda), Sendai, Japan

\end{slide}

%%%%%%%%%%%%%%%%%%%%%%%%%%%%%%%%%%%%%%%%%%%%%%%%%%%%%%%%%%%%%%%%%%%%%%%%%%%%%%%%%%%%%%%%%%%%%%%%%%%%%%%%%%%%%%%%%%%%%%%%%%%%%%%%%%%%%%%
%%%%%%%%%%%%%%%%%%%%%%%%%%%%%%%%%%%%%%%%%%%%%%%%%%%%%%%%%%%%%%%%%%%%%%%%%%%%%%%%%%%%%%%%%%%%%%%%%%%%%%%%%%%%%%%%%%%%%%%%%%%%%%%%%%%%%%%
%%%%%%%%%%%%%%%%%%%%%%%%%%%%%%%%%%%%%%%%%%%%%%%%%%%%%%%%%%%%%%%%%%%%%%%%%%%%%%%%%%%%%%%%%%%%%%%%%%%%%%%%%%%%%%%%%%%%%%%%%%%%%%%%%%%%%%%

\begin{slide}

(1996),
125--147 or Kyungpook Math. J. {\bf 40} (2000), 209--237 or
arXiv:math.NT/0602267.

[15] $ \textbf{J.-H. Yang}$, {\em The Method of Orbits for Real
Lie Groups,} Kyungpook Math. J. {\bf 42\, (2)} (2002), 199--272 or
arXiv:math.RT/0602056.

[16] $ \textbf{J.-H. Yang}$, \emph{Invariant metrics and Laplacians on Siegel-Jacobi spaces,}
Journal of Number Theory, {\bf 127} (2007), 83-102 or arXiv:math.
NT/0507215 v2.

[17] $\textbf{J.-H. Yang}$, \emph{The partial Cayley transform of
the Siegel-Jacobi disk ,} J. Korean Math. Soc. {\bf 45\,(3)} (2008), 781-794 or arXiv:math.NT\\ /0507216 v2.

[18] $ \textbf{J.-H. Yang}$, \emph{Invariant metrics and
Laplacians on Siegel-Jacobi disk,} arXiv:math. NT\\
/0507217 v2 or to appear in Chinese Annals of Math.

\end{slide}

%%%%%%%%%%%%%%%%%%%%%%%%%%%%%%%%%%%%%%%%%%%%%%%%%%%%%%%%%%%%%%%%%%%%%%%%%%%%%%%%%%%%%%%%%%%%%%%%%%%%%%%%%%%%%%%%%%%%%%%%%%%%%%%%%%%%%%%
%%%%%%%%%%%%%%%%%%%%%%%%%%%%%%%%%%%%%%%%%%%%%%%%%%%%%%%%%%%%%%%%%%%%%%%%%%%%%%%%%%%%%%%%%%%%%%%%%%%%%%%%%%%%%%%%%%%%%%%%%%%%%%%%%%%%%%%
%%%%%%%%%%%%%%%%%%%%%%%%%%%%%%%%%%%%%%%%%%%%%%%%%%%%%%%%%%%%%%%%%%%%%%%%%%%%%%%%%%%%%%%%%%%%%%%%%%%%%%%%%%%%%%%%%%%%%%%%%%%%%%%%%%%%%%%

\begin{slide}

[19] $ \textbf{J.-H. Yang}$, \emph{A note on a fundamental domain
for Siegel-Jacobi space}, Houston Journal of Mathematics, Vol.
{\bf 32}, No. 3 (2006), 701--712 or arXiv:math.NT/0507218.

[20] $ \textbf{J.-H. Yang}$, \emph{A note on invariant differential operators on Siegel-Jacobi spaces,} arXiv:math.\\
NT/0611388 v2 or revised version (2009).

%\end{slide}

%\newpage

\vskip 3cm

\begin{center}  \textbf{\Large Thank You Very Much !!!}\end{center}

\begin{center}
{\Large $ \textbf{$\bigstar\ \ \bigstar\ \ \bigstar\ \ \bigstar\ \ \bigstar
\ \ \bigstar\ \ \bigstar\ \ \bigstar$}$}
\end{center}

\begin{center}
{\Large $ \textbf{$\bigstar\ \ \bigstar\ \ \bigstar\ \ \bigstar\ \ \bigstar
\ \ \bigstar\ \ \bigstar\ \ \bigstar$}$}
\end{center}

\begin{center}
{\Large $ \textbf{$\bigstar\ \ \bigstar\ \ \bigstar\ \ \bigstar\ \ \bigstar
\ \ \bigstar\ \ \bigstar\ \ \bigstar$}$}
\end{center}

\begin{center}
{\Large $ \textbf{$\bigstar\ \ \bigstar\ \ \bigstar\ \ \bigstar\ \ \bigstar
\ \ \bigstar\ \ \bigstar\ \ \bigstar$}$}
\end{center}

\end{slide}

\end{document}